\documentclass[11pt,a4paper]{amsart}

\usepackage{amssymb}
\usepackage{amsfonts}
\usepackage{amsmath, amsthm}
\usepackage{fullpage}
\usepackage[all,cmtip]{xy}

\usepackage{multirow}
\usepackage{pgf}
\usepackage{yfonts}
\usepackage[colorlinks=true,linkcolor=blue,citecolor=red]{hyperref}

\usepackage{marginnote}

\theoremstyle{plain}
\newtheorem{theorem}{Theorem}[section]
\newtheorem{lemma}[theorem]{Lemma}
\newtheorem{proposition}[theorem]{Proposition}
\newtheorem{corollary}[theorem]{Corollary}
\theoremstyle{definition}
\newtheorem{definition}[theorem]{Definition}
\newtheorem{remark}[theorem]{Remark}

\newcommand{\N}{\mathbb{N}}
\newcommand{\Z}{\mathbb{Z}}
\newcommand{\Q}{\mathbb{Q}}
\newcommand{\R}{\mathbb{R}}

\newcommand{\K}{\mathbb{K}}

\newcommand{\calC}{\mathcal{C}}
\newcommand{\calD}{\mathcal{D}}

\newcommand{\calL}{\mathcal{L}}

\newcommand{\calK}{\mathcal{K}}

\newcommand{\Pj}{\mathbb{P}}
\newcommand{\Proj}{Proj}

\newcommand{\lieg}{\mathfrak{g}}
\newcommand{\lieb}{\mathfrak{b}}
\newcommand{\lieh}{\mathfrak{h}}

\DeclareMathOperator{\rlex}{rlex}
\DeclareMathOperator{\leqt}{\leq^t}

\DeclareMathOperator{\leqtrlex}{\leq^t_{\rlex}}

\DeclareMathOperator{\supp}{supp}
\DeclareMathOperator{\LS}{{\mathsf{LS}}}
\DeclareMathOperator{\SL}{{\mathsf{SL}}}
\DeclareMathOperator{\lieSL}{{\mathfrak{sl}}}

\newcommand{\pzero}{\textswab{0}}
\newcommand{\pone}{\textswab{1}}
\newcommand{\bond}{\mathfrak{b}}
\DeclareMathOperator{\gr}{gr}
\DeclareMathOperator{\NO}{\Delta}

\DeclareMathOperator{\conv}{conv}

\newcommand{\vleq}{\preceq}
\newcommand{\vlneq}{\prec}
\newcommand{\vgeq}{\succeq}
\newcommand{\vgneq}{\succ}
\newcommand{\lsleq}{\trianglelefteq}
\newcommand{\lslneq}{\vartriangleleft}

\newcommand{\lsgneq}{\vartriangleright}

\newcommand{\myref}[2]{\hyperref[#2]{#1 \ref{#2}}}

\title{LS algebras, valuations and Schubert varieties}

\author{Rocco Chiriv\`\i}
\address{Dipartimento di Matematica e Fisica ``Ennio De Giorgi'', Universit\`a del Salento, Lecce, Italy}
\email{rocco.chirivi@unisalento.it}

\author{Xin Fang}
\address{Mathematisches Institut, Universit\"at zu K\"oln, 50931, Cologne, Germany}
\email{xfang@math.uni-koeln.de}

\author{Peter Littelmann}
\address{Mathematisches Institut, Universit\"at zu K\"oln, 50931, Cologne, Germany}
\email{peter.littelmann@math.uni-koeln.de}


\begin{document}

\begin{abstract} In this paper, we propose an algebraic approach via Lakshmibai-Seshadri (LS) algebras to establish a link between standard monomial theories, Newton-Okounkov bodies and valuations. This is applied to Schubert varieties, where this approach is compatible with the one using Seshadri stratifications in \cite{CFL}, showing that LS paths encode vanishing multiplicities with respect to the web of Schubert varieties.
\end{abstract}

\maketitle

\section{Introduction}

\subsection{The case of Grassmannian varieties} In \cite{fangLittelmann}, the second and third authors studied the links between standard monomial theory, semi-toric degenerations of Grassmann varieties and the theory of Newton-Okounkov bodies. For this simple but already non-trivial case, they showed how the semi-toric degeneration of a Grassmann variety, compatible with all its Schubert subvarieties, may be defined in terms of a quasi-valuation on the field of rational functions of the variety.

Let $\K$ be an algebraically closed field and let $G_{d,n}$ be the Grassmann variety of $d$--dimensional subspaces of $\K^n$. Consider the set $I(d,n)$ of subsets of size $d$ in $\{1,2,\ldots,n\}$ and write $\tau \in I(d,n)$ as an increasing sequence $\tau = i_1 i_2 \cdots i_d$. We can partially order $I(d,n)$ by declaring $i_1 i_2 \cdots i_d \leq j_1 j_2 \cdots j_d$ if and only if $i_h \leq j_h$ for $h = 1,2, \ldots, d$. The Schubert varieties in $G_{d,n}$ are indexed by $I(d,n)$ and the Schubert variety $X(\sigma)$ is a codimension one subvariety of $X(\tau)$ if and only if $\tau$ covers $\sigma$ in $I(d,n)$.

Note that the set $I(d,n)$ indexes also the Pl\"ucker coordinates $p_\tau$ for the embedding $G_{d,n} \longrightarrow \Pj(\Lambda^d\K^n)$. With these generators the coordinate ring $A$ of the cone over $G_{d,n}$ in this embedding is a Hodge algebra (see \cite{dep}). This means that the set of monomials $p_{\tau_1}p_{\tau_2}\cdots p_{\tau_r}$ with $\tau_1\leq\tau_2\leq\cdots\leq\tau_r$, called standard, is a basis of $A$ as a $\K$--vector space and that each non-standard monomial $p_{\tau_1}p_{\tau_2}\cdots p_{\tau_r}$ is a linear combination of standard monomials $p_{\sigma_1}p_{\sigma_2}\cdots p_{\sigma_r}$ with $\tau_1\tau_2\cdots\tau_r$ lexicographically less than $\sigma_1\sigma_2\cdots\sigma_r$. These relations are called the straightening relations of $A$.

Let us write functions $I(d,n)\longrightarrow\Q$, i.e. elements of $\Q^{I(d,n)}$, as $\Q$--linear combinations of elements of $I(d,n)$. In \cite{fangLittelmann} the second and third author defined a valuation $\nu_\calC: A\setminus\{0\} \longrightarrow \Q^\calC\subseteq\Q^{I(d,n)}$ on $A$ associated to a maximal chain $\calC$ in $I(d,n)$ such that $\nu_\calC(p_\tau) = \tau$ for any $\tau\in\calC$. Next they defined a quasi-valuation $\nu$ by taking the minimum of the valuations $\nu_\calC$, as $\calC$ runs over the maximal chains of $I(d,n)$. This quasi-valuation $\nu$ induces a semi-toric degeneration of $G_{d,n}$ to a union of projective spaces $\Pj^N$, $N=\dim G_{d,n}$, one for each maximal chain $\calC$; moreover this degeneration is compatible with the Schubert varieties. The Newton-Okounkov body of a valuation $\nu_\calC$ is unimodular to the Gelfand-Tsetlin polytope and the standard monomials with support in $\calC$ define a simplex embedded in this polytope.

These results are obtained by exploiting the lattice structure of $I(d,n)$ and proving that certain standard monomials do appear in the straightening relations.

\subsection{From Grassmannians to Schubert varieties} It is natural to look for a generalization of the aforementioned results to partial flag varieties and related Schubert varieties. Let $G$ be a semisimple algebraic group and let $B\subseteq P$ be a Borel and, respectively, a parabolic subgroup of $G$. Let $\lambda$ be a dominant weight with stabilizer $P$, let $V_\lambda$ be the irreducible $G$ representation with highest weight $\lambda$ and consider the embedding of the partial flag variety $G/P$ in $\Pj(V_\lambda)$.

The Schubert varieties in $G/P$ are indexed by the set $W^\lambda$ of minimal representatives of the Weyl group $W$ of $G$ modulo the stabilizer of $\lambda$; recall that the Bruhat order of $W^\lambda$ is such that the Schubert variety $X(\sigma)$ is a codimension one subvariety of $X(\tau)$ if and only if $\tau$ covers $\sigma$ for this order. Now fix such a variety $X(\tau)$, with $\tau\in W^\lambda$ and consider the subposet $W^\lambda_\tau$ of the elements $\sigma$ less than or equal to $\tau$. So the Schubert subvarieties of $X(\tau)$ are parametrized by $W^\lambda_\tau$; hence this set is the analogue of $I(d,n)$ for the Grassmannian $G_{d,n}$ as far as the containment of Schubert varieties is concerned.

The analogue of the Pl\"ucker coordinates for the Schubert variety $X(\tau)$ is considerably more involved. We need the language of LS (Lakshmibai-Seshadri) paths as introduced by the third author in \cite{L1}, \cite{L2}, \cite{L3} and \cite{L4}. First we have to define certain positive integers $\bond(\sigma,\eta)$, called bonds, colouring the covering relations $\sigma < \eta$ in $W^\lambda_\tau$. An LS path is a function $\pi\in\Q^{W^\lambda_\tau}_{\geq 0}$ with totally ordered support and satisfying certain integral conditions with respect to the bonds along the support (see \myref{Section}{section_combinatoricsLSpaths}).

In \cite{L5} the third author, using quantum groups, defined certain generators $p_\pi$ of the homogeneous coordinate ring $A$ of (the cone over the embedding of) $X(\tau)$ in $\Pj(V_\lambda)$ associated to LS paths $\pi$. These generators $p_\pi$ take the place of the Pl\"ucker coordinates in this setting. We get again that the standard monomials in the $p_\pi$'s form a basis of $A$ as a vector space and each non-standard monomial $p_\pi p_{\pi'}$ is a linear combination of standard monomials $p_\eta p_{\eta'}$ with $\pi + \pi'$ less than $\eta + \eta'$ for a certain order on $\Q^{W^\lambda_\tau}$. These straightening relations have been proved in \cite{L5} and in \cite{LLM} by the third author and Lakshmibai, Magyar and the third author, respectively.

The combinatorial-algebraic content of this construction has been abstracted in the theory of LS algebras by the first author in \cite{chiriviLS} and \cite{chiriviLS2}. For details about LS algebras see \myref{Section}{section_LSAlgebra} and for the application to Schubert varieties see \myref{Section}{section_application} in the present paper.

Note that the case of Grassmannians corresponds to $G=\mathsf{SL}_n(\K)$, $\lambda$ the $d$--th fundamental weight and $\tau$ the longest element in $W^\lambda$; the poset $I(d,n)$ is isomorphic to $W^\lambda_\tau$, the bonds are all equal to $1$, so the LS paths (of degree $1$) are just the elements of $I(d,n)$ and the sections $p_\sigma$, $\sigma\in I(d,n)$ are, up to non-zero scalars, the Pl\"ucker coordinates.

\subsection{The geometric approach: Seshadri stratifications} In \cite{CFL} the authors gave a geometric construction of the quasi-valuation, described in the case of Grassmannians, for \emph{any} embedded projective variety that is non-singular in codimension one. In the case of Grassmannians studied in \cite{fangLittelmann} the Hodge algebra structure via the Pl\"ucker coordinates played a central role; however, in this new geometric approach the focus was shifted to the web of subvarieties (Schubert varieties for the Grassmannians).

Let $X$ be an embedded projective variety, let $(S,\leq)$ be a graded poset and let $(X_\sigma)_{\sigma\in S}$ be a collection of non-singular in codimension one subvarieties of $X$ indexed by $S$. Assume that $X_\sigma$ is a codimension one subvariety of $X_\tau$ for each $\sigma$ covered by $\tau$ in $S$ and that $S$ has a unique maximal element $\pone$ with $X_\pone = X$. Suppose moreover that we have a collection $(p_\sigma)_{\sigma\in S}$ of homogeneous functions such that: (1) $p_\sigma$ vanishes on $X_\tau$ if $\sigma\not\leq\tau$, (2) the set of zeros of $p_\tau$ in $X_\tau$ is set-theoretically the union of all the $X_{\sigma}$ such that $\sigma$ is covered by $\tau$. We call these data a \emph{Seshadri stratification} for $X$.

In \cite{CFL} a valuation $\nu_\calC:A\setminus\{0\} \longrightarrow \Q^\calC\subseteq\Q^S$ on the coordinate ring $A$ of the cone over $X$ is defined for each maximal chain $\calC$ in $S$ for a suitable order on $\Q^\calC$. This valuation is defined by looking at the order of vanishing of functions on $X_\tau$ along the divisor $X_\sigma$ for each covering step $\sigma < \tau$ in the chain $\calC$; these vanishing orders are normalized by the vanishing orders of the extremal functions $p_\tau$. A key property of $\nu_\calC$ is that $\nu_\calC(p_\tau) = \tau$ for each $\tau\in\calC$ (here we are using again the convention of writing elements of $\Q^S$ as linear rational combination of element of $S$). Moreover $\nu_\calC(f) > 0$ in $\Q^\calC$ for any $f\in A$ of positive degree.

Next a quasi-valuation $\nu$ is defined as the minimum, with respect to a certain total order on $\Q^S$, of the valuations $\nu_\calC$, with $\calC$ running over the set of all maximal chains of $S$. This construction is similar to that used in the theory of Newton-Okounkov bodies (see \cite{okounkov96}, \cite{kavehKhovanskii} and \cite{lazarsfeldMustata}) but instead of using a flag of subvarieties and uniformizers the quasi-valuation is defined in terms of the web of subvarieties $(X_\sigma)_{\sigma\in S}$ and the extremal functions $p_\sigma$, $\sigma\in S$.

Surprisingly, an application of Bertini's theorem guarantees that any embedded projective variety, which is non-singular in codimension one, admits a Seshadri stratification. 

As an example of application of this theory, in \cite{CFL}, \cite{CFL2} the case of Schubert varieties is studied. In particular, the quasi-valuation $\nu$ induces a semi-toric flat degeneration of a Schubert variety, compatible with all Schubert subvarieties, and the Newton-Okounkov body of $\nu$ is a geometric realization of the simplicial complex of chains in the Bruhat order.

\subsection{An algebraic approach via LS algebras} In the present paper we want to construct valuations and a quasi-valuation for an LS algebra similar to what seen before for Hodge algebras of Grassmannians and for projective varieties with a Seshadri stratification. The main combinatorial tool for this construction is the order requirement in the straightening relations. These relations should replace the geometric content of the Seshadri stratification. This construction is done in \myref{Section}{section_valuationLS}.

Let $(S,\leq)$ be a finite graded poset with a unique minimal element $\pzero$ and a unique maximal element $\pone$, let $\bond(\sigma,\tau)$ be a positive integer for each covering $\sigma < \tau$ in $S$ and let $\LS$ be the set of LS paths for the poset with bonds $(S,\leq,\bond)$ (see \myref{Section}{section_combinatoricsLSpaths} for details). Suppose that $A$ is an LS algebra over this poset with bonds.

In order to be able to construct the valuations and the quasi-valuation we need certain hypothesis on $A$. We say that $A$ has an effective system of weights if there exists a grading of $A$ by a free abelian group such that, for each $\sigma\in S$ and $r > 0$, the monomial $\sigma^r$ is the unique monomial of its degree. Moreover we say that $A$ is of flag type if: (1) it has an effective system of weights and (2) certain natural quotients of $A$ are domains. This requirement about quotients is geometrically clear: certain subvarieties of $\Proj(A)$ are irreducible, so we have a kind of stratification of $\Proj(A)$. (See \myref{Section}{section_LSAlgebra} for details.)

Now fix a maximal chain $\calC:\pzero = \sigma_0 < \sigma_1 < \cdots < \sigma_N = \pone$ in $S$. Denote by $\vleq$ the reverse lexicographic order on $\Q^\calC$. The first step is the construction of a valuation $\nu_\calC : A\setminus\{0\} \longrightarrow \Q^\calC$ such that $\nu_\calC(\tau) = \tau$ for each $\tau\in\calC$. In \myref{Theorem}{theorem_positivityForFlagType} we prove that such a valuation $\nu_\calC$ exists if $A$ is of flag type; moreover $\nu_\calC$ is a positive valuation: $\nu_\calC(f) \vgneq 0$ for any non-invertible $f\in A\setminus\{0\}$.

In \myref{Theorem}{theorem_estimate}, for an LS algebra with an effective system of weight, we prove an estimate for \emph{any} positive valuation $\nu_\calC$ such that $\nu_\calC(\tau) = \tau$ for each element $\tau$ in the fixed maximal chain $\calC$.

The final step is the definition of the quasi-valuation $\nu$ as the minimum (with respect to a suitable extension of $\vleq$ to $\Q^S$) of the valuations $\nu_\calC$'s with these properties, one for each maximal chain $\calC$ in $S$. The estimate in \myref{Theorem}{theorem_estimate} is crucial in proving that $\nu$ has the key property: $\nu(\pi) = \pi$ for every LS path $\pi$ (see \myref{Theorem}{theorem_quasiValuation}); hence $\nu$ has a strong link with the standard monomial theory of $A$ encoded in the LS structure.

In particular, at this point, we are finally ready to see the consequences of this construction in \myref{Section}{section_consequences}. We show that the quasi-valuation $\nu$ induces a flat degeneration of the LS algebra $A$ to the discrete LS algebra over $(S,\leq,\bond)$ and that the Newton-Okounkov body of $A$ with respect to $\nu$ is a geometric realization of the simplicial complex of chains of $(S,\leq)$.

\subsection{Application to Schubert varieties} In \myref{Section}{section_application} we apply this algebraic approach to Schubert varieties since the coordinate rings of the cone over the above defined embeddings are LS algebras of flag type: indeed (1) such rings have effective systems of weights induced by the torus action and (2) the Schubert subvarieties are irreducible. So we get a geometric interpretation of LS paths as image of the quasi-valuation $\nu$, a semi-toric degeneration result for Schubert varieties and the Newton-Okounkov body with respect to $\nu$ is the simplicial complex of chains in $W^\lambda_\tau$.

In \cite{kaveh} Kaveh defines valuations attached to certain maximal chains in $W^\lambda_\tau$ (those corresponding to simple roots) using the geometry of Bott-Samelson varieties. For minuscule $\lambda$ these valuations are equal to the valuation $\nu_\calC$'s defined in the present paper. But already for the adjoint representation of $\SL_3$ the valuations of Kaveh are different from the $\nu_\calC$'s.

Finally note that \myref{Theorem}{theorem_quasiValuation} can be applied also to the quasi-valuation defined for Schubert varieties via the Seshadri stratification in \cite{CFL}, so we get that $\nu(p_\pi) = \pi$ for each LS path; in this way we see that the geometric approach via the Seshadri stratification is coherent with the standard monomial theory for Schubert varieties. This last consequence is also proved in \cite{CFL2} using the construction of the sections $p_\pi$ via quantum groups in a more direct manner.

Note that we \emph{do not} prove that the valuation $\nu_\calC$ defined algebraically using the structure of LS algebra for the coordinate ring of a Schubert variety is equal to the valuation $\nu_\calC$ defined using the Seshadri stratification of the Schubert variety in \cite{CFL}. As far as our examples shows this is true and it would be of some interest to prove it in general.

\vskip 0.5 cm

{\bf Acknowledgements.} The first author would like to thank Andrea Maffei for countless useful discussions.

\section{Combinatorics of LS Paths}\label{section_combinatoricsLSpaths}

This Section introduces one of the main objects of the paper, the LS paths. The definition, given in \cite{chiriviLS}, has been modelled on that given by the third author in \cite{L2}. We study certain orders on LS paths, they will appear in the straightening relations of LS algebras. Then we introduce the order complex of a partially ordered set with its integral structure.

\subsection{Poset with bonds} Let $(S,\leq)$ be a finite partially ordered set, a \emph{poset} for short, with a unique minimal element $\pzero$ and a unique maximal element $\pone$; assume further that the poset is graded. Let $\sigma\in S$ and consider a maximal chain $\pzero=\sigma_0<\sigma_1<\cdots<\sigma_\ell=\sigma$ from $\pzero$ to $\sigma$ in $S$; being the poset graded, $\ell$ is independent of the chain, it is called the \emph{length} $\ell(\sigma)$ of $\sigma$. The length $\ell(\pone)$ of the maximal element of $S$ is called the length of $S$ and will always be denoted by $N$ in the sequel.

The element $\tau$ \emph{covers} the element $\sigma$ in $S$ if (1) $\sigma<\tau$ and (2) if $\sigma\leq\eta\leq\tau$ then either $\eta=\tau$ or $\eta=\sigma$.

For an element $\tau\in S$, we define $S_{\geq\tau} = \{\sigma\in S\,|\,\sigma\geq\tau\}$, $S_{>\tau} = \{\sigma\in S\,|\,\sigma > \tau\}$, $S_{\leq\tau} = \{\sigma\in S\,|\,\sigma \leq \tau\}$ and so on.

We say that an element $\tau$ is $\leq$--maximal with a certain property $P$ if $P(\tau)$ is true and $P(\sigma)$ is false for any $\sigma > \tau$.

\begin{definition}\label{definition_bonds} A \emph{set of bonds} on $S$ is a map $\bond$ from pairs $(\sigma,\tau)$, with $\tau$ covering $\sigma$, to positive integers such that: given two maximal chains
\[
\sigma = \eta_{i,1}<\eta_{i,2}<\cdots<\eta_{i,r} = \tau,\quad i = 1,2
\]
from $\sigma$ to $\tau$ in $S$ we have
\[
\gcd_{1\leq j\leq r-1}\bond(\eta_{1,j},\eta_{1,j+1}) \,\,= \gcd_{1\leq j\leq r-1}\bond(\eta_{2,j},\eta_{2,j+1}).
\]
\end{definition}

The map $\bond$ can be extended to comparable pairs $(\sigma,\tau)$ as the greatest common divisor on a maximal chain from $\sigma$ to $\tau$. For an element $\sigma\in S$, we denote by $M_\sigma$ the least common multiple of all bonds $\bond(\eta,\tau)$, with $\tau$ covering $\eta$, and either $\eta=\sigma$ or $\tau=\sigma$. Note that, by the way we have defined the bonds on comparable pairs, $M_\sigma$ is also the least common multiple of all bonds $\bond(\eta,\tau)$, with $\eta<\tau$ and either $\eta=\sigma$ or $\tau=\sigma$.

Fix an element $\tau\in S$. It is clear that the restriction $\bond_{|S_{\leq\tau}}$ to $S_{\leq\tau}$ of a set of bond $\bond$ is a set of bonds for the poset $(S_{\leq\tau},\leq_{|S_{\leq\tau}})$ with minimal element $\pzero$ and maximal element $\tau$. The same is true for $S_{\geq\tau}$.

\subsection{LS paths} Let $(S,\leq,\bond)$ be a poset with bonds.

\begin{definition}\label{definition_lspath} A function $\pi:S\longrightarrow\Q$ is an \emph{LS path} of \emph{degree} $\deg\pi = r\in\N$ if the following conditions hold:
\begin{itemize}
\item[(1)] $\pi(\sigma)\geq0$ for any $\sigma\in S$,
\item[(2)] $\supp\pi$, defined as the set of $\sigma\in S$ such that $\pi(\sigma) \neq 0$, is a totally ordered subset of $S$,
\item[(3)] if $\supp\pi=\{\sigma_1<\sigma_2<\cdots<\sigma_n\}$ then
\[
\bond(\sigma_j,\sigma_{j+1})\sum_{i=1}^j\pi(\sigma_i)\in\N\quad\textrm{for all }j = 1, \ldots, n-1,
\] and
\[
\sum_{i=1}^n\pi(\sigma_i)=r.
\]
\end{itemize}
\end{definition}

It follows at once from the definition that $M_\sigma\pi(\sigma)$ is an integer for any LS path $\pi$ and any $\sigma\in S$. Note that in (3) we may equivalently suppose that $\supp\pi\subseteq\{\sigma_1<\sigma_2<\cdots<\sigma_n\}$. Further, if $\pi$ is a real valued function on $S$ satisfying the above condition (3), then $\pi(\sigma)\in\Q$ for any $\sigma\in S$; so the conditions (1), (2) and (3) characterize the LS paths of degree $r$ in the set of real valued functions on $S$.

Fix an element $\tau \in S$. The LS paths of $(S_{\leq\tau}, \leq_{|S_{\leq\tau}}, \bond_{|S_{\leq\tau}})$ can be clearly identified with the LS paths $\pi$ of $(S,\leq,\bond)$ such that $\max\supp\pi\leq\tau$.

We write elements of $\Q^S$, in particular LS paths, as linear rational combinations of elements of $S$ by identifying $\sigma\in S$ with the function
\[
S\ni\eta\longmapsto\left\{\begin{array}{ll}1 & \textrm{if } \eta = \sigma\\0 & \textrm{otherwise}\end{array}\right.\in\Q
\]
so that, for any such function $\pi$ we have:
\[
\pi = \sum_{\sigma\in S}\pi(\sigma)\sigma.
\]
If $\sigma$ is an element of $S$ and $r\in\N$, then $r\sigma$ is called the \emph{extremal} LS path of degree $r$ and support $\{\sigma\}$. The \emph{width} $w(\pi)$ of an LS path $\pi$ is defined as $w(\pi) = \ell(\max\supp\pi) - \ell(\min\supp\pi) + 1$; in particular the extremal LS paths are those of width $1$.

We denote by $\LS_r$ the set of all LS paths of degree $r$ and by $\LS$ the set of all LS paths. Given a chain $\calC$ in $S$, we denote by $\LS(\calC)$ the set of LS paths having support contained in $\calC$, accordingly $\LS_r(\calC)$ is the set of LS paths of degree $r$ in $\LS(\calC)$. We say that the LS paths $\pi_1,\pi_2,\ldots,\pi_r$ have \emph{comparable supports} if there exists a chain in $S$ containing $\supp\pi_1\cup\cdots\cup\supp\pi_r$; i.e., there exists a chain $\calC$ such that $\pi_1,\pi_2,\ldots,\pi_r\in\LS(\calC)$. In such a case the sum $\pi_1+\pi_2+\cdots+\pi_r$ is an LS path of degree $\deg\pi_1 + \cdots + \deg\pi_r$.
\begin{proposition}[Proposition~3 in \cite{chiriviLS}]
If $\pi$ is an LS path of degree $r$ then there exist LS paths $\pi_1,\pi_2,\cdots,\pi_r$ of degree $1$, such that: $\max\supp\pi_h\leq\min\supp\pi_{h+1}$, for all $h=1,2,\ldots,r-1$, and $\pi=\pi_1 + \pi_2 + \cdots + \pi_r$ as functions on $S$.
\end{proposition}

\begin{proposition}[Remark after Proposition~1.2 in \cite{chiriviLS}]\label{proposition_normal_lattice} Let $\calC$ be a chain in $S$, let $L_\calC$ be the lattice in $\R^\calC$ generated by $\LS(\calC)$. Then $\LS(\calC)$ is the set of non-negative functions in $L_\calC$.
\end{proposition}

A (formal) monomial $\pi_1\pi_2\cdots\pi_r$ of LS paths is \emph{standard} if $\max\supp\pi_h\leq\min\supp\pi_{h+1}$ for $h=1,2,\ldots,r-1$. If the LS paths $\pi_1,\pi_2,\ldots,\pi_r$ have comparable supports then $\pi=\pi_1+\cdots+\pi_r$ is an LS path and, by \myref{Proposition}{proposition_normal_lattice}, there exist LS paths $\pi_{0,1},\pi_{0,2},\dots,\pi_{0,s}$ of degree $1$ such that (1) $\pi=\pi_{0,1}+\cdots+\pi_{0,s}$ and (2) the monomial $\pi_{0,1}\cdots\pi_{0,s}$, called the \emph{canonical form} of $\pi_1\pi_2\cdots\pi_r$ and of $\pi$, is standard.

We state here some simple results about LS paths, they are immediate consequences of the definitions.

\begin{lemma}\label{lemma_splitLSpath} Let $\pi\in\LS_1$ with $\supp\pi=\{\sigma_1<\ldots<\sigma_n\}$, $a_i=\pi(\sigma_i)$ for $i=1,\ldots,n$ and $n\geq 2$. Then $\pi'=(1-a_1)\sigma_1 + a_1\sigma_2$ and $\pi''=(a_1+a_2)\sigma_2 + \pi_{|\{\sigma_3,\ldots,\sigma_n\}}$ are LS paths of degree $1$. Moreover $\sigma_1\pi''$ is the canonical form of the non-standard monomial $\pi'\pi$.
\end{lemma}

\begin{lemma}\label{lemma_w2LSpath} If $a\sigma + b\tau$ is an LS path of degree $1$, then $b\sigma + a\tau$ is also an LS path of degree $1$. If moreover $a\geq 1/2$ then $(2a-1)\sigma + 2b\tau$ is an LS paths of degree $1$.
\end{lemma}

\subsection{Orders}\label{subsection_orders} Let $\leqt$ be a total order on $S$ refining the partial order $\leq$. The set $\Q^S$ of all functions $\pi:S\longrightarrow\Q$ can be totally ordered using the reverse lexicographic order: we define
\[
\pi\leqtrlex\pi'
\]
if and only if either $\pi=\pi'$ or, denoting by $\sigma$ the $\leqt$--maximal element with $\pi(\sigma)\neq\pi'(\sigma)$, we have $\pi(\sigma)<\pi'(\sigma)$. We define also a partial order $\lsleq$ on $\Q^S$ by requiring $\pi\lsleq\pi'$ if $\pi\leqtrlex\pi'$ for any total order $\leqt$ refining $\leq$ on $S$.

We begin by seeing an equivalent form of the order $\lsleq$ we have just introduced; this equivalence will be used in the sequel often without explicit mention.

\begin{lemma}\label{lemma_orderEquivalence} Let $\pi,\pi'\in\Q^S$. The following are equivalent:
\begin{itemize}
	\item[(i)] $\pi\lslneq\pi'$,
	\item[(ii)] if $\tau$ is $\leq$--maximal such that $\pi(\tau)\neq\pi'(\tau)$ then $\pi(\tau) < \pi'(\tau)$.
\end{itemize}
\end{lemma}
\begin{proof} We prove that (i) implies (ii). So suppose that $\pi\lslneq\pi'$ and let $\tau$ be $\leq$--maximal with $\pi(\tau)\neq\pi'(\tau)$.

We define an order $\leq^t$ on $S$ as follows: (1) choose an arbitrary refinement of $\leq$ to a total order on $S_{\geq\tau}$, (2) choose an arbitrary refinement of $\leq$ to a total order on $S\setminus S_{\geq\tau}$ and (3) declare $\sigma<^t\tau$ for any $\sigma\in S\setminus S_{\geq\tau}$. It is clear that $\leq^t$ is a total order on $S$ refining $\leq$.

Note that $\tau$ is the $\leq^t$--minimal element of $S_{\geq\tau}$. Moreover $\pi(\sigma) = \pi'(\sigma)$ for all $\sigma\in S_{\geq\tau}\setminus\{\tau\}$ by our assumption on $\tau$. So $\tau$ is $\leq^t$--maximal such that $\pi(\tau) \neq \pi'(\tau)$, hence $\pi(\tau) < \pi'(\tau)$ by $\pi\lslneq\pi'$.

Now we prove that (ii) implies (i). Let $\leq^t$ be a total order on $S$ refining $\leq$ and let $\tau$ be $\leq^t$--maximal such that $\pi(\tau)\neq\pi'(\tau)$. Since $\leq^t$ refines $\leq$, the element $\tau$ is also $\leq$--maximal with the property that $\pi(\tau)\neq\pi'(\tau)$, hence $\pi(\tau) < \pi'(\tau)$ by (ii). This shows that $\pi\lslneq\pi'$.
\end{proof}

Now we state some properties of the order $\lsleq$ that will be used in the study of the properties of valuations related to LS paths.
\begin{lemma}\label{lemma_orderCompatibleWithSum}
The order $\lsleq$ is compatible with the addition of functions, i.e. if $\pi,\pi',\eta\in\Q^S$ and $\pi\lsleq\pi'$ then $\pi + \eta\lsleq\pi' + \eta$.
\end{lemma}
\begin{proof} This is clear by definition.
\end{proof}
\begin{lemma}\label{lemma_orderMaxSupp} Let $\pi,\pi'$ be LS paths. If $\pi\lsleq\pi'$ then $\max\supp\pi\leq\max\supp\pi'$.
\end{lemma}
\begin{proof}
Let $\tau = \max\supp\pi$ and $\sigma > \tau$. If $\pi'(\sigma) > 0$ then $\sigma\in\supp\pi'$ and $\max\supp\pi'\geq\sigma > \tau$ and our claim is proved. So we may assume that $\pi'_{|S_{>\tau}} = 0$.

Now, if we had $\pi'(\tau) = 0$ then $\tau$ would be $\leq$--maximal such that $\pi(\tau)\neq\pi'(\tau)$ and we have $\pi(\tau) > 0 = \pi'(\tau)$. This is impossible: \myref{Lemma}{lemma_orderEquivalence} forces $\pi(\tau) < \pi'(\tau)$ since $\pi \lslneq \pi'$.

Hence we must have $\pi'(\tau) > 0$ and, using again that $\pi'_{|S_{>\tau}} = 0$, we have $\max\supp\pi' = \tau$ and our claim is proved.
\end{proof}

\begin{lemma}\label{lemma_orderOnSum} Let $\pi_1,\pi_2,\pi_1',\pi_2'$ be LS paths of degree $1$. If $\pi_1 + \pi_2\lsleq\pi_1' + \pi_2'$ and $\pi_1'\pi_2'$ is a standard monomial then: $\pi_2\lsleq\pi_2'$ and, if $\pi_2=\pi_2'$ then $\pi_1\lsleq\pi_1'$.
\end{lemma}
\begin{proof} We prove the first claim. If $\pi_2=\pi_2'$ then the claim is true; so suppose $\pi_2\neq\pi_2'$ and let $\tau$ be $\leq$--maximal such that $\pi_2(\tau) \neq \pi_2'(\tau)$. Suppose by contradiction that $\pi_2(\tau) > \pi_2'(\tau)$. We would have
\[
1\geq\sum_{\sigma\in S_{\geq\tau}}\pi_2(\sigma) > \sum_{\sigma\in S_{\geq\tau}}\pi_2'(\sigma),
\]
so $\supp\pi_2'\not\subseteq S_{\geq\tau}$. Hence, being $\pi_1'\pi_2'$ standard, the support of $\pi_1'$ does not intersect $S_{\geq\tau}$. For each $\sigma\in S_{>\tau}$ we find
\[
(\pi_1 + \pi_2)(\sigma) \geq \pi_2(\sigma) = \pi_2'(\sigma) = (\pi_1' + \pi_2')(\sigma)
\]
and since we have also
\[
(\pi_1 + \pi_2)(\tau) \geq \pi_2(\tau) > \pi_2'(\tau) = (\pi_1' + \pi_2')(\tau),
\]
we conclude $\pi_1 + \pi_2 \lsgneq \pi_1' + \pi_2'$ that is impossible. This finishes the proof of the first claim.

The second part of the claim follows by \myref{Lemma}{lemma_orderCompatibleWithSum}.
\end{proof}

\subsection{Order complex}\label{subsection_orderComplex}

The LS paths over the poset with bond $(S,\leq,\bond)$ are maps $S\longrightarrow\Q$. We can consider them as elements of $\R^S$ whose basis is the set of extremal LS paths $\sigma$, with $\sigma\in S$, of degree $1$. In the following definition we use the Euclidean topology of $\R^S$.
\begin{definition}\label{definition_orderComplex} The \emph{order complex} of $S$ is defined as
\[
\Delta(S)=\overline{\Big\{\frac{\pi}{r}\,|\,\pi\in\LS_r,\,r > 0\Big\}}\subseteq\R^S
\]
and, given a chain $\calC$ in $S$, the corresponding \emph{simplex} is
\[
\Delta(\calC)=\overline{\Big\{\frac{\pi}{r}\,|\,\pi\in\LS_r(\calC),\,r > 0\Big\}}\subseteq\R^S.
\]
\end{definition}
These names are justified by the following proposition.
\begin{proposition}\label{proposition_chainSimplex}
For a chain $\calC$ the set $\Delta(\calC)$ is the $(|\calC|-1)$--dimensional simplex with set of vertices $\calC$; in particular, $\Delta(\calC)$ is a convex set. Moreover $\Delta(S)$ is a simplicial complex with faces the simplices $\Delta(\calC)$, $\calC$ a chain of $S$.
\end{proposition}
\begin{proof} We fix a chain $\calC$ and we prove that $D\doteq\{\pi/r\,|\,\pi\in\LS_r(\calC),\,r > 0\}$ is closed under rational convex linear combinations; from this it clearly follows that $\Delta(\calC) = \overline{D}$ is a convex set.

Let $\pi\in\LS_r(\calC)$, $\pi'\in\LS_{r'}(\calC)$ and let
\[
\eta = \frac{a}{a+b} \cdot \frac{\pi}{r} + \frac{b}{a+b} \cdot \frac{\pi'}{r'}
\]
be a convex combination with $a,b\in\N$. Let $\overline{\pi}=ar'\pi$ and $\overline{\pi}'=br\pi'$; these are LS paths of degree $arr'$ and $brr'$, respectively; hence $\overline{\pi}+\overline{\pi}'$ is an LS path of degree $(a+b)rr'$. We have $\eta = (\overline{\pi}+\overline{\pi}')/(a+b)rr'\in D$.

The rest of the proof is very neat. The elements $\sigma$, $\sigma\in\calC$, belong to $\Delta(\calC)$, hence their convex hull, i.e. the $(|\calC|-1)$--simplex having them as vertices, is contained in $\Delta(\calC)$. On the other hand, if $\pi\in\LS_r(\calC)$, then $\pi/r=\sum_{\sigma\in\calC}\pi(\sigma)\sigma/r$ is a convex combination of the elements $\sigma$'s, $\sigma\in\calC$; hence $\pi/r$ is an element of their convex hull. This finishes the proof of the first statement.

The second statement is a consequence of the first one since it is clear that: (1) $\Delta(S)$ is the union of the $\Delta(\calC)$'s, $\calC$ a chain of $S$, and (2) $\Delta(\calC_1)\cap\Delta(\calC_2) = \Delta(\calC_1\cap\calC_2)$ for any pair of chains $\calC_1$, $\calC_2$.
\end{proof}

\begin{corollary}\label{corollary_orderGeometric} The order complex $\Delta(S)$ is a geometric realization in $\R^S$ of the simplicial complex of chains in $(S,\leq)$.
\end{corollary}
\begin{proof} This follows at once by the previous proposition and the definition of the simplicial complex of chains (see, for example, \cite{bjorner}).
\end{proof}

The next step is the introduction of a finer structure on $\Delta(S)$, in order to reflect the bonds. Following Raika Dehy's paper \cite{dehy}, we recall the definition of an \emph{integral structure} on a simplicial complex $\calK\subseteq\R^n$: it is a collection $\calK_r$, $r=1,2,\ldots$, of subsets of $\calK$ such that for each simplex $K$ of $\calK$ there exists an affine embedding $i_K:K\longrightarrow\R^k$, with $k=\dim K$, such that (a) $i_K(K)$ has vertices with integral coordinates and (b) $\calK_r\cap K = i_K^{-1}(\frac{1}{r}\Z^k\cap i_K(K))$, for each $r=1,2,\ldots$. We want to show that we have a natural integral structure on the simplicial complex $\Delta(S)$ using the bonds of $S$. The proof of the following proposition is an adaptation of the analogous result in \cite{dehy}.
\begin{proposition}\label{proposition_integralStructureViaBonds}
The collection of subsets
\[
\Delta_r(S)=\Big\{\frac{\pi}{r}\,|\,\pi\in\LS_r\Big\},\quad r=1,2,\ldots
\]
is an integral structure for $\Delta(S)$. Moreover $\Delta(S) = \overline{\bigcup_{r\geq 1}\Delta_r(S)}$.
\end{proposition}
\begin{proof} Fix the simplex $K=\Delta(\calC)$ of the chain $\calC:\sigma_0<\sigma_1<\cdots<\sigma_k$ of $S$, and note that $\Delta_r(S)\cap K$ is the set of all functions $\pi/r$ with $\pi\in\LS_r(\calC)$. We define the embedding $i_K$ as the linear extension of the map
\[
\sigma_i \longmapsto v_i \doteq \sum_{j=1}^i\bond(\sigma_{j-1},\sigma_j)e_j, \quad i = 0,1,\ldots,k,
\]
where $e_1,\ldots,e_k$ is the canonical basis of $\R^k$. It is then clear that the vertices of $i_K(K)$, i.e. the images of $\sigma_0,\ldots,\sigma_k$, have integral coordinates. We need to show that, for a function $\varphi:\calC \longrightarrow \R$, we have $i_K(\varphi/r)\in\frac{1}{r}\Z^k\cap i_K(K)$ if and only if $\varphi\in\LS_r(\calC)$.

Let $\varphi=\sum_{i=0}^ka_i\sigma_i$ and $\bond(\sigma_{j-1},\sigma_j)=n_j$ for short, then
\[
\begin{array}{rcl}
i_K(\displaystyle\frac{\varphi}{r}) & = & \displaystyle\frac{1}{r}\sum\limits_{i=0}^k a_i v_i\\
 & = & \displaystyle\frac{1}{r}\sum\limits_{j=1}^k n_j \left(\sum\limits_{i=j}^k a_i\right) e_j.
\end{array}
\]
But, being $v_1,\ldots, v_k$ linearly independent, the first equality implies that $i_K(\varphi/r)\in i_K(K)$ if and only if $a_i\geq0$, for any $i=0,1,\ldots,k$, and $\sum_{i=0}^ra_i=r$. So we find $i_K(\varphi/r)\in\frac{1}{r}\Z^k\cap i_K(K)$ if and only if $a_i\geq0$, for any $i$, $\sum_{i=0}^r a_i=r$ and, by the second equality above, $n_j\sum_{i=j}^k a_i = n_jr - n_j\sum_{i=0}^{j-1}a_i\in\Z$ for $j=1,\ldots,k$. This last condition is clearly equivalent to $n_j\sum_{i=0}^{j-1}a_i\in\Z$, for $j=1,\ldots,k$. But these conditions are those in the definition of LS paths; this finishes the proof.
\end{proof}

\section{LS Algebras}\label{section_LSAlgebra}

In this section we introduce LS algebras and their straightening relations. We define a discrete LS algebra as an LS algebra with the simplest straightening relations and we recall how certain quotients of an LS algebra are still LS algebras. Further we define what is an LS algebra of flag type, this kind of algebra is going to play a crucial role in the sequel. For more details and results about LS algebras see \cite{chiriviLS} and \cite{chiriviLS2}.

\smallskip

Let $\K$ be a field, let $A$ be a commutative $\K$--algebra, fix an injection
\[
\LS_1\,\ni\,\pi\,\longmapsto\,\pi\,\in\,A
\]
and extend it to $\LS$ using the canonical form for higher degree LS paths: if $\pi=\pi_1 + \pi_2 + \cdots + \pi_r$ is the canonical form of $\pi$, then we map
\[
\pi\,\longmapsto\, \pi_1\pi_2\cdots\pi_r \in\,A.
\]

\begin{definition}\label{definition_lsalgebra} The algebra $A$ is an \emph{LS algebra} over the poset with bonds $(S,\leq,\bond)$ if
\begin{enumerate}
\item[(LS1)] the set of LS paths is a basis of $A$ as a $\K$--vector space and the degree of LS paths induces a grading for $A$,
\item[(LS2)] if $\pi_1\pi_2$ is a non-standard monomial of degree $2$ and
\[
\pi_1\pi_2 = \sum_j c_j\pi_{j,1}\pi_{j,2}
\]
is the unique relation, called a \emph{straightening relation}, expressing $\pi_1\pi_2$ as a $\K$--linear combination of standard monomials, as guaranteed by (LS1), then $\pi_1 + \pi_2\lsleq\pi_{j,1} + \pi_{j,2}$, for any $j$ such that $c_j\neq0$,
\item[(LS3)] if $\pi_1\pi_2$ is a non-standard monomial of degree $2$ and $\pi_1$, $\pi_2$ have comparable supports then the canonical form of $\pi_1\pi_2$ appears with coefficient $1$ in the straightening relation for $\pi_1\pi_2$.
\end{enumerate}
The LS algebra $A$ is \emph{discrete} if the straightening relations are
\[
\pi_1\pi_2 = \left\{
\begin{array}{ll}
\pi_{0,1}\pi_{0,2} & \textrm{if } \pi_{0,1}\pi_{0,2} \textrm{ is the canonical form of }\pi_1\pi_2,\\[0.5em]
0 & \textrm{if } \pi_1 \textrm{ and } \pi_2 \textrm{ have non-comparable supports.}\\
\end{array}
\right.
\]
\end{definition}

In the sequel, if not otherwise stated, we write a straightening relation as above with implicitly assuming that $c_j\neq0$ for all $j$.

We remark that if $\pi_1,\pi_2$ have comparable supports, as in (LS3), and $\pi_{0,1}\pi_{0,2}$ is the canonical form of the non-standard monomial $\pi_1\pi_2$, then $\pi_{0,1} + \pi_{0,2}=\pi_1 + \pi_2$, hence $\pi_{0,1}\pi_{0,2}$ is the $\lsleq$--minimal element appearing in the straightening relation for $\pi_1\pi_2$. The above definition may be slightly generalized allowing the coefficient of $\pi_{0,1}\pi_{0,2}$ to be any non-zero scalar; all the results of this paper hold also for this generality. (Compare with the definition of special LS algebra in \cite{chiriviLS} and Theorem~11.1 in \cite{CFL}.)

We call the set of elements $\pi\in A$, for $\pi\in\LS_1$, a \emph{path basis} for the LS algebra $A$. Note that an algebra may have various path basis.

If all bonds equal $1$ then an LS algebra is a special type of a Hodge algebra (see \cite{dep}).

In the following Lemma we state some consequences of the straightening relations for an LS algebra; this will be used in inductive arguments.

\begin{lemma}\label{lemma_inequalityInRelation} If $\pi_1'\pi_2'$ is a standard monomial appearing in the straightening relation for the non-standard monomial $\pi_1\pi_2$, then $\pi_1,\pi_2 \lslneq \pi_2'$.
\end{lemma}
\begin{proof} We may consider only $i=2$ by symmetry. The inequality $\pi_2\lsleq\pi_2'$ follows by \myref{Lemma}{lemma_orderOnSum} using the order requirement in (LS2). Suppose, by contradiction, that $\pi_2=\pi_2'$. Then, again by \myref{Lemma}{lemma_orderOnSum}, $\pi_1\lsleq\pi_1'$ and we derive $\max\supp\pi_1\leq\max\supp\pi_1'$ by \myref{Lemma}{lemma_orderMaxSupp}. Being $\pi_1'\pi_2'$ standard we have also $\max\supp\pi_1'\leq\min\supp\pi_2'=\min\supp\pi_2$, hence $\pi_1\pi_2$ is a standard monomial, contrary to the assumption.
\end{proof}

As proved in Proposition~8 in \cite{chiriviLS}, the straightening relations in (LS2) and (LS3) for degree $2$ monomials are equivalent to:
\begin{enumerate}
\item[(LS2')] if $\pi_1\pi_2\cdots\pi_r$ is a non-standard monomial of degree $r\geq2$ and
\[
\pi_1\pi_2\cdots\pi_r=\sum_j c_j\pi_{j,1}\pi_{j,2}\cdots\pi_{j,r}
\]
is the unique relation expressing $\pi_1\pi_2\cdots\pi_r$ as a $\K$--linear combination of standard monomials, as guaranteed by (LS1), then $\pi_1 + \pi_2 + \cdots + \pi_r \lsleq \pi_{j,1} + \pi_{j,2} + \cdots + \pi_{j,r}$ for any $j$ such that $c_j\neq0$,
\item[(LS3')] if $\pi_1\pi_2\cdots\pi_r$ is a non-standard monomial of degree $r\geq2$ and $\pi_1$, $\pi_2,\ldots,\pi_r$ have comparable supports then the canonical form of $\pi_1 + \pi_2 + \cdots + \pi_r$ appears with coefficient $1$ in the straightening relation for $\pi_1\pi_2\cdots\pi_r$.
\end{enumerate}

\begin{definition}\label{definition_weights}
Let $\Lambda$ be a free abelian group and let $\lambda:\LS\longrightarrow{\Lambda}$ be such that $\lambda(\pi_1 + \pi_2) = \lambda(\pi_1) + \lambda(\pi_2)$ for each pair of LS paths $\pi_1$, $\pi_2$. Using the injection $\LS\longrightarrow A$ we can define $\lambda$ for the path basis of the LS algebra $A$. We say that $\lambda$ is a \emph{system of weights} for $A$ if the straightening relations are $\lambda$--homogeneous.

The system of weights $\lambda$ is \emph{effective} if, for any $\sigma\in S$ and $r\in\N$, the monomial $\sigma^r$ is the unique LS path monomial of degree $r$ with weight $r\lambda(\sigma)$.
\end{definition}

\begin{lemma}\label{lemma_weights}
If $\lambda$ is a system of weights for $A$ then
\[
\lambda(\pi) = \sum_{\sigma\in S}\pi(\sigma)\lambda(\sigma),
\]
in particular $\lambda$ is completely determined by its values on the set of extremal paths $\sigma\in S$.
\end{lemma}
\begin{proof} Let $\supp\pi = \{\sigma_1, \ldots,\sigma_n\}$ and let $M$ be such that $M\pi(\sigma_i) \in \N$ for all $i = 1,2,\ldots,n$. Note that the monomial $\sigma_1^{M\pi(\sigma_1)}\cdots\sigma_n^{M\pi(\sigma_n)}$ appears in the straightening relation for the monomial $\pi^M$. Hence these two monomials have the same weight and the claim follows since $\Lambda$ is free.
\end{proof}

We now recall how certain quotients of an LS algebra are still LS algebras. Fix an element $\tau\in S$ and define $I_\tau$ as the ideal of the LS algebra $A$ generated by the LS paths $\pi$ such that $\max\supp\pi\not\leq\tau$.
\begin{proposition}[see Proposition 23 in \cite{chiriviLS}]\label{proposition_quotientLS} The quotient $A_\tau\doteq A/I_\tau$ is an LS algebra over the poset with bonds $(S_{\leq\tau}, \leq_{|S_\tau}, \bond_{|S_\tau})$ with injection $\LS(S_{\leq\tau})\ni\pi\longmapsto \pi + I_\tau \in A_\tau$.
\end{proposition}

The following definition is central for the whole paper.
\begin{definition}\label{definition_flag}
We say that an LS algebra $A$ is of \emph{flag type} if $A$ has an effective system of weights and $A_\tau$ is a domain for any $\tau$ in $S$.
\end{definition}

\begin{remark}\label{remark_stillFlag}
Note that $I_\pone=0$, hence an LS algebra of flag type is a domain. Note also that $\pi + I_\tau \longmapsto \lambda(\pi) \in \Lambda$ is an effective system of weights for $A_\tau$ if $\lambda$ is an effective system of weights for $A$. Hence $A_\tau$ is still of flag type for any $\tau\in S$ if $A$ is of flag type.
\end{remark}

\section{Valuations and Newton-Okounkov Bodies}\label{section_NewtonOkounkovBody}

We introduce quasi-valuations and the associated Newton-Okounkov bodies. We use a definition of Newton-Okounkov body inspired by \cite{okounkov03} since it is more adequate to our aims for a quasi-valuation. When this quasi-valuation is a valuation on a graded domain, our definition coincides with the one in \cite{kavehKhovanskii}, \cite{lazarsfeldMustata} and \cite{okounkov96}. The definition of the levels of the Newton-Okounkov body is modelled on that of integral structure for a simplicial complex in \myref{Subsection}{subsection_orderComplex}. For more details about valuations see, for example, \cite{englerPrestel}.

\begin{definition}\label{definition_quasiValuation} A \emph{quasi-valuation} on an algebra $A$ over a field $\K$ with values in a totally ordered abelian group $V$ is a map $\nu:A\setminus\{0\}\longrightarrow V$ such that
\begin{itemize}
\item[(a)] $\nu(x+y)\geq\min\{\nu(x),\nu(y)\}$ for all $x,y\neq0$ such that $x+y\neq0$,
\item[(b)] $\nu(cx)=\nu(x)$ for all $x\neq0$ and $c\in\K^*$,
\item[(c)] $\nu(xy)\geq\nu(x)+\nu(y)$ for all $x,y\neq0$ such that $xy\neq0$.
\end{itemize}
\noindent If moreover in the property (c) we have $\nu(xy) = \nu(x)+\nu(y)$, for all $x,y,xy\neq0$, then we say that $\nu$ is a \emph{valuation}.

\noindent A quasi-valuation is \emph{positive} if $\nu(x) > 0$ for every non-zero and non-invertible $x\in A$.
\end{definition}

If the algebra $A$ is a domain and $\nu:A\setminus\{0\}\longrightarrow V$ is a valuation, then we can extend $\nu$ to a valuation on the quotient field $Q(A)$ of $A$ by declaring: $\nu(x/y) = \nu(x) - \nu(y)$ for all $x,y\in A\setminus\{0\}$.

An immediate consequence of the above properties (a) and (b) is the following.
\begin{lemma}\label{lemma_quasiValuationMin}
If $x_1,x_2,\ldots,x_n\neq0$ and $\nu(x_1)<\nu(x_2),\ldots,\nu(x_n)$ then $\nu(x_1+\cdots+x_n)=\nu(x_1)$.
\end{lemma}

In the sequel we will construct certain quasi-valuations as in the following proposition; its easy proof is omitted.
\begin{proposition}\label{proposition_quasiValuationAsMin}
If $\nu_1,\nu_2,\ldots,\nu_n:A\setminus\{0\}\longrightarrow V$ are valuations on an algebra $A$, then $A\setminus\{0\}\ni x\longmapsto\nu(x)=\min\{\nu_1(x),\nu_2(x),\ldots,\nu_n(x)\}\in V$ is a quasi-valuation on $A$.
\end{proposition}

\begin{definition}\label{definition_nobody} Suppose that $V$ is a $\Z$--submodule of a real vector space of finite dimension. Given a quasi-valuation $\nu:A\setminus\{0\}\longrightarrow V$ of a graded algebra $A$, we define the \emph{Newton-Okounkov body} of $A$ with respect to $\nu$ as
\[
\NO_\nu(A) \doteq \overline{\left\{\frac{1}{r}\nu(x)\,|\,x\in A_r\setminus\{0\},\,r=1,2,\ldots\right\}}
\]
where the closure is with respect to the euclidean topology of the real vector space. For $r\geq 1$ the subset
\[
\left\{\frac{1}{r}\nu(x)\,|\,x\in A_r\setminus\{0\}\right\}\subseteq\NO_\nu(A)
\]
is called the $r$--th \emph{level} of $\NO_\nu(A)$.
\end{definition}
The notion of level will play an important role in the connection to combinatorics of LS algebras in later Sections. Note that the Newton-Okounkov body is the closure of the union of all the levels.

In the following Proposition we show that $\NO_\nu(A)$ is a convex set for a valuation on a graded domain $A$; this recovers the usual definition of Newton-Okounkov body (see for example \cite{kavehKhovanskii}, \cite{lazarsfeldMustata} and \cite{okounkov96}) up to a translation.
\begin{proposition}\label{proposition_noQuasiNo}
Let $\nu$ be a valuation for the graded domain $A$. Then $\NO_\nu(A)$ is a convex set. If moreover $x_1\in A$ is any element of degree $1$, then
\[
\NO_\nu(A) - \nu(x_1) = \overline{\conv\left\{\frac{1}{r}\nu(\frac{x}{x_1^r})\,|\,x\in A_r\setminus\{0\},\,r=1,2,\ldots\right\}}.
\]
\end{proposition}
\begin{proof} We just need to show that $\NO_\nu(A)$ is a convex set. This follows once we prove that $C = \{\nu(x)/r\,|\,x\in A_r\setminus\{0\},\, r=1,2,\ldots\}$ is closed under rational convex combination.

So let $u = \nu(x)/r$ and $v = \nu(y)/s$ be two elements of $C$ and let $w = a u/c + b v/c$, with $a+b=c$, $a,b\in\N$, be a rational convex combination of them. Set $e=as$ and $f=br$ and consider the non-zero element $z = x^ey^f$ of degree $crs$. We find at once $\nu(z)/crs = (e\nu(x) + f\nu(y))/crs = a\nu(x)/cr + b\nu(y)/cs = au/c + bv/c$. (Compare with the proof of \myref{Proposition}{proposition_chainSimplex}).
\end{proof}

We add some words of explanation about our definition of $\NO_\nu$; one should compare this definition with the second formula in Section 3.10 of Okounkov's paper \cite{okounkov03}. We want to use the Newton-Okounkov body also for quasi-valuations. The main problem with the right hand side formula in the previous Proposition is the convexification operator: it is not necessary for valuations on graded domains but it totally spoils the Newton-Okounkov body for quasi-valuations.

\section{Valuations for LS Algebras}\label{section_valuationLS}

In this core section we first prove that an LS algebra of flag type has a positive valuation with prescribed valued along a fixed chain. Next we prove an estimate for such valuations that we crucially use later in the proof of the properties of the quasi-valuation defined as the minimum of the valuations along the chains; such properties show the links of the quasi-valuation with the standard monomial theory.

\subsection{Existence of positive valuations}\label{subsection_existencePositive}

Let $(S,\leq,\bond)$ be a poset with bonds, recall that we denote by $N$ the length of $S$. In this and the next subsection we fix once and for all a maximal chain $\calC:\pzero=\sigma_0<\sigma_1<\cdots<\sigma_N=\pone$ in $S$.

Our aim is to define a valuation on $A$ with values in $\Q^\calC$; so we need to define a total order on $\Q^\calC$. We use the reverse lexicographic order induced by the total order of $\calC$ and denote it by $\vleq$. Note that, when restricted to $\Q^\calC$, the order $\vleq$ coincides with $\lsleq$.

\begin{theorem}\label{theorem_positivityForFlagType} Suppose that $A$ is an LS algebra of flag type, then there exists a positive valuation $\nu$ on $A$ with values in $\Q^\calC$ such that $\nu(\sigma_h) = \sigma_h$ for all $h = 0, \ldots, N$.
\end{theorem}
\begin{proof} We denote by $\lambda$ an effective system of weights for $A$. We proceed by induction on the length $N=\ell(\pone)$ of $S$. If $N=0$, then $S=\{\sigma_0\}$ and $A$ is the polynomial ring $\K[\sigma_0]$, so the assignment $\nu(\sigma_0) = \sigma_0$ may be extended to the lowest term valuation which is clearly positive and has values in $\Q^\calC$.

Now suppose $N>0$, for short let $\sigma=\sigma_{N-1}$ be the second to last element in $\calC$ and $\tau=\sigma_N$ the last one. Further let $I_\sigma$ be the ideal of $A$ generated by the LS paths $\pi$ with $\max\supp\pi\not\leq\sigma$ and recall that $A_\sigma=A/I_\sigma$ is an LS algebra over the poset $S_{\leq\sigma}=\{\tau\in S\,|\,\tau\leq\sigma\}$ by \myref{Proposition}{proposition_quotientLS} and it is also of flag type by \myref{Remark}{remark_stillFlag}. Note that $\calC':\pzero=\sigma_0<\cdots<\sigma_{N-1}=\sigma$ is a maximal chain in $S_{\leq\sigma}$. Hence, by induction, there exists a positive valuation $\nu'$ on $A_\sigma$, with values in $\Q^{\calC'}\subseteq \Q^\calC$, such that $\nu'(\sigma_h) = \sigma_h$, for $h=0,1,\ldots,N-1$.

We want to extend $\nu'$ to a valuation $\nu$ on $A$. In order to construct $\nu$ we consider the multiplicative set $M = \{1,\sigma,\sigma^2,\ldots\}$ and we claim that $M^{-1}I_\sigma$ is a principal ideal of $M^{-1}A$. Indeed let $n=\bond(\sigma,\tau)$ be the last bond in the chain $\calC$ and let $\eta_1$ be the LS path
\[
\eta_1 =  (1 - \frac{1}{n})\sigma + \frac{1}{n}\tau,
\]
we want to show that $M^{-1}I_\sigma = M^{-1}A\eta_1$. Since $\max\supp\eta_1 = \tau$, it is clear that $\eta_1\in I_\sigma$. We use (reverse) induction on $\lsleq$, for $\pi\in\LS_1$ such that $\max\supp\pi\not\leq\sigma$, to show that $\pi\in M^{-1}A\eta_1$.

Consider the monomial $\sigma\pi$, we use two different strategies depending on whether this is a standard or a non-standard monomial. Suppose first that $\sigma\pi$ is standard, then $\sigma\leq\min\supp\pi$ since the other inequality $\max\supp\pi\leq\sigma$ is impossible. So
\[
\pi = \eta_k = (1-\frac{k}{n})\sigma + \frac{k}{n}\tau
\]
for some $1\leq k\leq n$. We prove that $\eta_k\in M^{-1}A\eta_1$, for any $k\geq 1$, by (reverse) induction on $k$; for $k=n$ we have $\eta_n=\tau$, hence this proves also the basic step for the induction on $\lsleq$ since $\tau$ is the $\lsleq$--maximal LS path.

We can clearly suppose $k\geq2$. Recall that by (LS2') there exists a straightening relation
\[
\eta_1^k = \sum_j c_j\pi_{j,1}\cdots\pi_{j,k}
\]
for the non-standard monomial $\eta_1^k$. The canonical form of $\eta_1^k$ is $(k-1)\sigma + \eta_k$ since $\sigma\leq\min\supp\eta_k$, so $(k-1)\sigma+\eta_k\lsleq\pi_{j,1} + \cdots + \pi_{j,k}$ for any $j$. This inequality implies $\eta_k\lsleq\pi_{j,k}$ by the definition of $\lsleq$. If this inequality is strict then $\pi_{j,k}\in M^{-1}A\eta_1$ by induction on $\lsleq$. On the other hand, if $\eta_k=\pi_{j,k}$, then, using the $\lambda$--homogeneity of the above straightening relation, we find that the standard monomial $\pi_{j,1}\cdots\pi_{j,k-1}$ of degree $k-1$ has weight $(k-1)\lambda(\sigma)$. But $\lambda$ is effective, hence $\pi_{j,1}\cdots\pi_{j,k-1}=\sigma^{k-1}$. Moreover, the monomial $\sigma^{k-1}\eta_k$ does appear with coefficient $1$ by (LS3'). So we have proved $\eta_1^k=\sigma^{k-1}\eta_k + x$, where $x\in M^{-1}A\eta_1$ is the sum of all the monomials in the straightening relation with $\pi_{j,k}\neq\eta_k$, and we conclude $\eta_k=(\eta_1^k - x)/\sigma^{k-1}\in M^{-1}A\eta_1$. This finishes the proof that $\pi\in M^{-1}A\eta_1$ if $\sigma\pi$ is standard.

Note that for $k=n$, the straightening relation has the simpler form $\eta_1^n=\sigma^{n-1}\tau$; we will need this later in the proof.

\smallskip

Now suppose that $\sigma\pi$ is a non-standard monomial, then in its straightening relation
\[
\sigma\pi = \sum_j c_j\pi_{j,1}\pi_{j,2}
\]
we have $\pi\lslneq\pi_{j,2}$ for any $j$. Indeed, if otherwise $\pi = \pi_{j,2}$, then $\lambda(\sigma\pi) = \lambda(\sigma) + \lambda(\pi) = \lambda(\pi_{j,1}) + \lambda(\pi)$ and, being $\lambda$ effective, $\pi_{j,1} = \sigma$ that is impossible since it would imply that $\sigma\pi = \pi_{j,1}\pi_{j,2}$ was standard.

Note that $\max\supp\pi\leq\max\supp\pi_{j,2}$ by \myref{Lemma}{lemma_orderMaxSupp}, so if $\max\supp\pi_{j,2}\leq\sigma$, then $\max\supp\pi\leq\sigma$ that is impossible. We find $\max\supp\pi_{j,2}\not\leq\sigma$ and $\pi_{j,2}\in I_\sigma$, hence by induction, $\pi_{j,2}\in M^{-1}A\eta_1$, for any $j$, and we have $\pi=(\sum_j c_j\pi_{j,1}\pi_{j,2})/\sigma\in M^{-1}A\eta_1$, concluding the proof of our claim that $M^{-1}I_\sigma = M^{-1}A\eta_1$.

At this point we know that $M^{-1}I_\sigma$ is prime, since $A_\sigma$ is a domain, and it is also a principal ideal; so the localization $M^{-1}A_{M^{-1}I_\sigma}$ is a discrete valuation ring with valuation
\[
M^{-1}A_{M^{-1}I_\sigma}\setminus\{0\} \ni x \stackrel{\overline{\nu}}{\longmapsto} r\in\N
\]
where $r$ is such that $x\in M^{-1}I_\sigma^r\setminus M^{-1}I_\sigma^{r+1}$. Note that if $x\in A\setminus\{0\}$ and $\overline{\nu}(x)=r$, then $x/\eta_1^r\not\in M^{-1}I_{\sigma}$, hence $x/\eta_1^r + I_\sigma$ is not zero in $Q(A_\sigma)$. We use this to define
\[
A\setminus\{0\}\ni x\stackrel{\nu}{\longmapsto}\nu'(x/\eta_1^r + I_\sigma) + r\left((1 - \frac{1}{n})\sigma + \frac{1}{n}\tau\right),
\]
where the right hand side belongs to $\Q^\calC$ since $\nu'(x/\eta_1^r + I_\sigma) \in \Q^{\calC'}$. Note that, for $x\not\in I_\sigma$, we have $r=0$, hence $\nu(x)=\nu'(x)$. Moreover $\overline{\nu}(\tau)=n$ since, as proved above, $\tau=\eta_1^n/\sigma^{n-1}$, thus 
\[
\nu(\tau) = \nu'(1/\sigma^{n-1}) + n\left((1 - \frac{1}{n})\sigma + \frac{1}{n}\tau\right) = \tau.
\]
So $\nu(\sigma_h) = \sigma_h$ for $h=0,1,\ldots,N$. We need to show that $\nu$ is a valuation, but this is a trivial check using that $\overline{\nu}$ and $\nu'$ are valuations.

Finally note that, for $x\in A\setminus\{0\}$, if $\overline{\nu}(x)>0$, then $\tau$ appears with a positive coefficient in $\nu(x)$, hence $\nu(x) \vgneq 0$; on the other hand if $\overline{\nu}(x)=0$, then $\nu(x)=\nu'(x) \vgneq 0$ if $x$ is not invertible by the inductive hypothesis. This proves that $\nu$ is a positive valuation.
\end{proof}

\subsection{Estimate of valuations for LS algebras}\label{subsection_estimate}

In order to estimate certain set of valuations in $\Q^\calC$ we use the following map
\[
\LS\, \ni\, \pi\, \stackrel{\nu_{0,\calC}}{\longmapsto}\,\sum_{\sigma\in S}\pi(\sigma)\sigma_{\ell(\sigma)}\,\in\,\Q^\calC.
\]
This is clearly an additive map: $\nu_{0,\calC}(\pi_1 + \pi_2) = \nu_{0,\calC}(\pi_1) + \nu_{0,\calC}(\pi_2)$. Moreover $\nu_{0,\calC}(\pi) = \pi$ if and only if $\supp\pi\subseteq\calC$.

The following Lemma will be needed in the proof of \myref{Theorem}{theorem_estimate}, we state it here for clarity.

\begin{lemma}\label{lemma_nu0Inequality} Let $\pi_1\pi_2$ be a standard monomial of degree $2$ and let $\pi_2'\in\LS_1$ be such that $\pi_2\lslneq\pi_2'$. Then $\nu_{0,\calC}(\pi_1 + \pi_2) \vlneq \nu_{0,\calC}(\pi_2')$.
\end{lemma}
\begin{proof}
Let $\leqt$ be a total order on $S$ refining $\leq$ and let $\tau$ be $\leqt$--maximal in $\supp\pi'_2$ such that $\pi_2(\tau)<\pi_2'(\tau)$. If $\tau \leqt \min\supp\pi_2$ then we would have
\[
1 = \sum_{\eta\in\supp\pi_2}\pi_2(\eta) = \pi_2(\tau) + \sum_{\eta\in\supp\pi_2,\,\eta > \tau}\pi_2(\eta) < \pi_2'(\tau) + \sum_{\eta\in\supp\pi_2,\,\eta > \tau}\pi_2'(\eta)\leq 1;
\]
so $\tau>^t\min\supp\pi_2$. Using $\max\supp\pi_1\leq\min\supp\pi_2<^t\tau$, we find $(\pi_1 + \pi_2)(\tau)=\pi_2(\tau)<\pi_2'(\tau)$, thus we conclude
\[
\begin{array}{rcl}
\nu_{0,\calC}(\pi_2') & \vgeq & \sum_{\eta > \tau} \pi_2'(\eta)\sigma_{\ell(\eta)} + \pi_2'(\tau)\sigma_{\ell(\tau)}\\[0.5em]
 & = & \sum_{\eta > \tau} \pi_2(\eta)\sigma_{\ell(\eta)} + \pi_2'(\tau)\sigma_{\ell(\tau)}\\[0.5em]
 & \vgneq & \sum_{\eta\in S} (\pi_1(\eta) + \pi_2(\eta))\sigma_{\ell(\eta)}\\[0.5em]
 & = & \nu_{0,\calC}(\pi_1 + \pi_2).
\end{array}
\]
\end{proof}

Given an LS path $\pi$ such that $\supp\pi\not\subseteq\calC$ let
\[
h_\calC(\pi) = \min\Big\{h\,|\,0\leq h\leq N \textrm{ and }\sigma_h > \max(\supp\pi\setminus\calC)\Big\},
\]
while, for an LS path with support in $\calC$, let $h_\calC(\pi)=-1$. We define also $\sigma_{-1}=0$, $M_{\sigma_{-1}} = 1$ and, for $-1\leq h\leq N$, $\sigma'_h = \sigma_h / M_{\sigma_h}$, $V_{h,\calC}=\Q^{\{\sigma_0, \ldots, \sigma_h\}} \subseteq \Q^\calC$. Note that $V_{-1,\calC} = \Q^\varnothing = \{0\}$. In the sequel we will simply write $h(\pi)$ and $V_h$ when the chain $\calC$ is clear from the context and this does not raise any ambiguity.

Another Lemma we will need is the following.
\begin{lemma}\label{lemma_impossibleCase} Let $\pi$ be an LS path of degree $1$ such that $\supp\pi\not\subseteq\calC$ and let $\pi_1\pi_2$ be a standard monomial of degree $2$. Let $h = h(\pi)$ and suppose that: (1) $\sigma_{h-1} + \pi \lslneq \pi_1 + \pi_2$ and (2) $\pi \lslneq \pi_2$. Then there exists $\tau\in S$ with $\ell(\tau) \geq h$ and $\pi(\tau) < \pi_2(\tau)$. Moreover such $\tau$ can be chosen to be $\leq$--maximal with the property $\pi(\tau) \neq \pi_2(\tau)$.
\end{lemma}
\begin{proof} We proceed by contradiction to prove the first claim of the Lemma. Let $S_{\geq h} = \{\tau\in S\,|\, \ell(\tau) \geq h\}$ and suppose that $\pi(\tau) \geq \pi_2(\tau)$ for all $\tau\in S_{\geq h}$.

We begin by ruling out the possibility that $\pi(\tau) > \pi_2(\tau)$ for some element $\tau\in S_{\geq h}$. Indeed, by choosing a $\leq$--maximal $\tau$ with this property we would have $\pi \lsgneq \pi_2$ while $\pi \lslneq \pi_2$ by hypothesis. So, in what follow, we assume $\pi_{|S_{\geq h}} = \pi_{2|S_{\geq h}}$.

Let $\epsilon = \max(\supp\pi\setminus\calC)$ and note that $\ell(\epsilon) < h$ and $\sigma_{h-1}$ and $\epsilon$ are non-comparable by the definition of $h$.

Let $\tau$ be $\leq$--maximal such that $(\sigma_{h-1} + \pi)(\tau) \neq (\pi_1 + \pi_2)(\tau)$. Such $\tau$ exists since $\sigma_{h-1} + \pi \neq \pi_1 + \pi_2$ by hypothesis. In particular $(\sigma_{h-1} + \pi)(\tau) < (\pi_1 + \pi_2)(\tau)$ by \myref{Lemma}{lemma_orderEquivalence} since $\sigma_{h-1} + \pi \lslneq \pi_1 + \pi_2$.

If $\ell(\tau) \geq h$ then $(\sigma_{h-1} + \pi)(\tau) = \pi(\tau) = \pi_2(\tau)$, hence $\pi_1(\tau) > 0$. So $\tau\in\supp\pi_1$, hence $\tau \leq \max\supp\pi_1 \leq \min\supp\pi_2$ being $\pi_1\pi_2$ standard, and we find that $\supp\pi_2\subseteq S_{\geq h}$. But this is impossible since we would have
\[
1 = \sum_{\sigma\in S_{\geq h}}\pi_2(\sigma) = \sum_{\sigma\in S_{\geq h}}\pi(\sigma) < 1
\]
where in the last inequality we have used $\pi(\epsilon) > 0$ and $\ell(\epsilon) < h$.

So we have $\ell(\tau) < h$. Now we want to show that $\pi_2(\sigma_{h-1}) > 0$.

If $(\sigma_{h-1} + \pi)(\sigma_{h-1}) \neq (\pi_1 + \pi_2)(\sigma_{h-1})$ then, using $\ell(\tau) < h$, we find that $\sigma_{h-1}$ is $\leq$--maximal with this property, hence $(\sigma_{h-1} + \pi) (\sigma_{h-1}) < (\pi_1 + \pi_2)(\sigma_{h-1})$ by \myref{Lemma}{lemma_orderEquivalence} since $\sigma_{h-1} + \pi \lslneq \pi_1 + \pi_2$. So, in any case, $1 \leq (\sigma_{h-1} + \pi)(\sigma_{h-1}) \leq (\pi_1 + \pi_2)(\sigma_{h-1})$.

But then if we had $\pi_2(\sigma_{h-1}) = 0$ then $\pi_1(\sigma_{h-1}) > 0$, hence $\sigma_{h-1} < \min\supp\pi_2$ so that $\supp\pi_2\subseteq S_{\geq h}$ which, as we have already seen, is impossible. This proves that $\pi_2(\sigma_{h-1}) > 0$ as claimed.

Now $\pi_2(\epsilon) = 0$ since $\pi_2(\sigma_{h-1}) > 0$ and $\epsilon$ and $\sigma_{h-1}$ are non-comparable while $\supp\pi_2$ is a totally ordered subset. In particular $\pi(\epsilon) \neq \pi_2(\epsilon)$. Let $\eta\geq\epsilon$ be $\leq$--maximal with this property.

It is clear that $\ell(\eta) \leq h-1$ since $\pi_{|S_{\geq h}} = \pi_{2|S_{\geq h}}$ by hypothesis. If $\eta\in\supp\pi_2$ we would have $\epsilon \leq \eta \leq \sigma_{h-1}$ since $\pi_2(\sigma_{h-1}) > 0$ and $\supp\pi_2$ is totally ordered. So $\pi_2(\eta) = 0$ since $\epsilon$ and $\sigma_{h-1}$ are non-comparable.

This is a contradiction: $\eta$ is $\leq$--maximal such that $\pi(\eta) \neq \pi_2(\eta)$ and $\pi(\eta) > \pi_2(\eta)$ whereas $\pi \lslneq \pi_2$ and by \myref{Lemma}{lemma_orderEquivalence} we should have $\pi(\eta) < \pi_2(\eta)$.

We have proved the first claim of the Lemma. The second one follows at once since $\pi \lslneq \pi_2$.
\end{proof}

In the following Theorem we prove an estimate that will play a key role in the next \myref{Section}{subsection_quasiValuation} to prove that the quasi-valuation we will study has a close link to the standard monomial theory.
\begin{theorem}\label{theorem_estimate} Let $A$ be a domain and an LS algebra with an effective system of weights and suppose that $\nu$ is a $\Q^\calC$--valued positive valuation such that $\nu(\sigma_h) = \sigma_h$, for $h=0,1,\ldots,N$. Then for any LS path $\pi$ we have
\begin{itemize}
\item[(a)] if $\supp\pi\subseteq\calC$ then $\nu(\pi) = \nu_{0,\calC}(\pi) = \pi$,
\item[(b)] if $\supp\pi\not\subseteq\calC$ then there exists $v'\in V_{h(\pi)-1}$ such that $\nu(\pi) \vgeq \nu_{0,\calC}(\pi) + \sigma'_{h(\pi)} + v'$.
\end{itemize}
In particular, $\nu(\pi) \vgeq \nu_{0,\calC}(\pi)$ for any LS path $\pi$; the equality holds if and only if $\supp\pi\subseteq\calC$.
\end{theorem}
\begin{proof} Note first that if $\supp\pi\subseteq\calC$ and we have proved $\nu(\pi)=\nu_{0,\calC}(\pi)$, then we have also $\nu(\pi) \vgeq \nu_{0,\calC}(\pi) + \sigma'_{h(\pi)} + v'$ with $v'=0$, indeed $h(\pi)=-1$, $\sigma'_{-1}=0$ and $V_{-1} = \{0\}$. So (b) is true also when $\supp\pi\subseteq\calC$; moreover in this case $v'=0$ is the unique possible choice.

We proceed in several steps for proving claims (a) and (b).

\vskip 0.125cm

\noindent{\bf$\bigstar$} Step 1. We reduce the claims to the case of LS paths of degree $1$. So assume (a) and (b) are true for $\LS_1$ and let $\pi$ be an LS path of degree $r$, $\pi=\pi_1\pi_2\cdots\pi_r$ its canonical form. Then $\nu(\pi)=\nu(\pi_1\cdots\pi_r)=\nu(\pi_1) + \cdots + \nu(\pi_r)$; on the other hand, $\pi=\pi_1 + \pi_2 + \cdots + \pi_r$ as functions on $S$, so $\nu_{0,\calC}(\pi)=\nu_{0,\calC}(\pi_1) + \nu_{0,\calC}(\pi_2) + \cdots + \nu_{0,\calC}(\pi_r)$.

Now, if $\supp\pi\subseteq\calC$, then $\supp\pi_i\subseteq\calC$ for each $i=1,\ldots,r$, hence $\nu(\pi)=\nu(\pi_1) + \cdots + \nu(\pi_r) = \nu_{0,\calC}(\pi_1) + \cdots + \nu_{0,\calC}(\pi_r)=\nu_{0,\calC}(\pi)$. On the other hand, if $\supp\pi\not\subseteq\calC$ let $j$ be maximal such that $\supp\pi_j\not\in\calC$. We know that  $\nu(\pi_i) \vgeq \nu_{0,\calC}(\pi_i) + \sigma'_{h(\pi_i)} + v'_i$, with $v'_i\in V_{h(\pi_i) - 1}$, for $i = 1,\ldots,r$. Moreover $h(\pi)=h(\pi_j)$ by definition, so we conclude $\nu(\pi) \vgeq \nu_{0,\calC}(\pi) + \sigma'_{h(\pi)} + v'$ with $v'=v_1'+\cdots+v'_j\in V_{h(\pi) - 1}$ since $v'_{j+1} = v'_{j+2} = \cdots = v'_r = 0$ by $\supp\pi_{j+1}\cup\cdots\cup\supp\pi_r\subseteq\calC$.

\smallskip

In the following steps we use reverse induction on $\lsleq$ to prove claims (a) and (b) for LS paths of degree 1. Note that $\pi=\sigma_N$ is the unique $\lsleq$--maximal LS path of degree $1$ and clearly $\nu(\pi) = \sigma_N = \nu_{0,\calC}(\pi)$.

\vskip 0.125cm

\noindent{\bf$\bigstar$} Step 2. In this step we prove claim (a) for LS path $\pi\in\LS_1$ with $\supp\pi\subseteq\calC$ and $w(\pi)\leq2$. If $w(\pi)=1$ the claim is clear. So let $\supp\pi=\{\sigma_h<\sigma_k\}$ and $\pi(\sigma_h)=a_h$, $\pi(\sigma_k)=a_k$, $a_h,a_k>0$ and $a_h + a_k = 1$.

\noindent{\bf$\star$} Step 2.1. Suppose that $a_h\geq 1/2$. By \myref{Lemma}{lemma_w2LSpath}, $\pi'=(2a_h-1)\sigma_h + 2a_k\sigma_k$ is an LS path of degree $1$ and $2\pi=\sigma_h + \pi'$, so $\sigma_h\pi'$ is the canonical form of the non-standard monomial $\pi^2$. Let
\[
\pi^2=\sigma_h\pi' + \sum_j c_j\pi_{j,1}\pi_{j,2}
\]
be the straightening relation for $\pi^2$, with $2\pi=\sigma_h + \pi' \lslneq \pi_{j,1} + \pi_{j,2}$ for any $j$. By \myref{Lemma}{lemma_inequalityInRelation}, $\pi\lslneq\pi'$ and $\pi\lslneq\pi_{j,2}$, for any $j$; so, by induction, $\nu(\pi')=\nu_{0,\calC}(\pi')$, since $\supp\pi'\subseteq\calC$, and $\nu(\pi_{j,2}) \vgeq \nu_{0,\calC}(\pi_{j,2})$. Moreover $\pi'\lsleq\pi_{j,2}$ by \myref{Lemma}{lemma_orderOnSum}. If $\pi'=\pi_{j,2}$ then, denoting by $\lambda$ the effective system of weights of $A$, we have $\lambda(\pi_{j,1}) + \lambda(\pi_{j,2}) = 2\lambda(\pi) = \lambda(\sigma_h) + \lambda(\pi')$, hence $\lambda(\pi_{j,1})=\lambda(\sigma_h)$ and so $\pi_{j,1}=\sigma_h$; this is clearly impossible since it would imply that $\sigma_h\pi' = \pi_{j,1}\pi_{j,2}$. So $\pi'\lslneq\pi_{j,2}$.

But then $\nu_{0,\calC}(\pi_{j,2}) \vgneq \nu_{0,\calC}(\sigma_h + \pi')$ by \myref{Lemma}{lemma_nu0Inequality} (applied to $\pi_1=\sigma_h$, $\pi_2=\pi'$, $\pi_2'=\pi_{j,2}$ with notation as in that Lemma) and, using that $\nu$ is positive in the first inequality below, we find
\[
\begin{array}{c}
\nu(\pi_{j,1}\pi_{j,2}) = \nu(\pi_{j,1}) + \nu(\pi_{j,2}) \vgeq \nu(\pi_{j,2}) \vgeq \nu_{0,\calC}(\pi_{j,2})\vgneq\\
\vgneq \nu_{0,\calC}(\sigma_h+\pi') = \nu_{0,\calC}(\sigma_h) + \nu_{0,\calC}(\pi') = \nu(\sigma_h) + \nu(\pi') = \nu(\sigma_h\pi').\\
\end{array}
\]
Hence we conclude $2\nu(\pi)=\nu(\sigma_h\pi')=\nu_{0,\calC}(\sigma_h\pi')=2\nu_{0,\calC}(\pi)$, where the first equality holds since, as we have just proved, $\sigma_h\pi'$ is the unique monomial appearing in the straightening relation with the lowest value of $\nu$.

\noindent{\bf$\star$} Step 2.2. Suppose that $a_h<1/2$ and consider $\pi'=a_k\sigma_h + a_h\sigma_k$, this is an LS path of degree $1$ by \myref{Lemma}{lemma_w2LSpath}. The monomial $\pi'\pi$ is non-standard, $\sigma_h\sigma_k$ is its canonical form and $\nu(\pi')=\nu_{0,\calC}(\pi')$ as proved in Step 2.1. In particular we have
\[
\pi'\pi = \sigma_h\sigma_k + \sum_j c_j\pi_{j,1}\pi_{j,2}.
\]
with $\sigma_h + \sigma_k \lslneq \pi_{j,1} + \pi_{j,2}$. Note that, by \myref{Lemma}{lemma_inequalityInRelation}, $\pi\lslneq\sigma_k$ and $\pi\lslneq\pi_{j,2}$, for any $j$. Moreover $\sigma_k\lsleq\pi_{j,2}$ by \myref{Lemma}{lemma_orderOnSum}, without equality since $\sigma_h$ is extremal and we have an effective system of weights. So, as before,
\[
\begin{array}{c}
\nu(\pi_{j,1}\pi_{j,2}) \vgeq \nu(\pi_{j,2}) \vgeq \nu_{0,\calC}(\pi_{j,2}) \vgneq\\
\vgneq \nu_{0,\calC}(\sigma_h\sigma_k) = \sigma_h + \sigma_k = \nu(\sigma_h\sigma_k).\\
\end{array}
\]
Hence $\nu(\pi'\pi) = \nu(\sigma_h\sigma_k) = \sigma_h + \sigma_k$. We conclude $\nu(\pi)=\nu(\sigma_h\sigma_k)-\nu(\pi')=\nu_{0,\calC}(\pi)$.

\vskip 0.125cm

\noindent{\bf$\bigstar$} Step 3. Now we consider the case of LS path $\pi$ of degree $1$ with $\supp\pi\subseteq\calC$ and $w(\pi)\geq3$ and prove (a) for $\pi$. Let $\supp\pi=\{\tau_1<\cdots<\tau_n\}$, with $\pi(\tau_i)=a_i$, for $i=1,\ldots,n$ and consider $\pi'=(1-a_1)\tau_1 + a_1\tau_2$, $\pi''=(a_1+a_2)\tau_2 + \pi_{|\{\tau_3,\ldots,\tau_n\}}$; by \myref{Lemma}{lemma_splitLSpath}, $\pi'$ and $\pi''$ are both LS paths of degree $1$ and $\tau_1\pi''$ is the canonical form of the non-standard monomial $\pi'\pi$. Let
\[
\pi'\pi = \tau_1\pi'' + \sum_j c_j\pi_{j,1}\pi_{j,2}
\]
be the straightening relation for $\pi'\pi$; note that $\tau_1 + \pi' \lslneq \pi_{j,1} + \pi_{j,2}$. As above, $\pi\lslneq\pi''$, also $\pi''\lsleq\pi_{j,2}$, for any $j$, and the equality is impossible since $\tau_1$ is an extremal path and $A$ has an effective system of weights. So, again as in the previous cases, we find $\nu(\pi'\pi) = \nu(\tau_1\pi'') = \nu_{0,\calC}(\tau_1\pi'')$ since $\pi'' \lsgneq \pi$. Now note that $\nu(\pi') = \nu_{0,\calC}(\pi')$, since $w(\pi')=2$ and $\supp\pi'\subseteq\calC$. We conclude $\nu(\pi) = \nu_{0,\calC}(\tau_1\pi'')-\nu_{0,\calC}(\pi') = \nu_{0,\calC}(\pi)$.

\vskip 0.125cm

\noindent{\bf$\bigstar$} Step 4. In this final step we assume that $\supp\pi\not\subseteq\calC$. Let $h=h(\pi)$. By definition $h\geq1$ and the monomial $\sigma_{h-1}\pi$ is non-standard. Moreover $\sigma_{h-1}$ and $\pi$ have non-comparable supports so in the straightening relation
\[
\sigma_{h-1}\pi = \sum_j c_j\pi_{j,1}\pi_{j,2}
\]
there is no canonical form but, being $A$ a domain, the right hand side is not zero. By (LS2) and \myref{Lemma}{lemma_inequalityInRelation}, $\pi\lslneq\pi_{j,2}$, for any $j$, hence, by induction, $\nu(\pi_{j,2}) \vgeq \nu_{0,\calC}(\pi_{j,2}) + \sigma'_{h(\pi_{j,2})} + v_j' \vgeq \nu_{0,\calC}(\pi_{j,2})$, with $v_j'\in V_{h(\pi_{j,2})-1}$.

Fix $j$ such that $\nu(\pi_{j,1}\pi_{j,2})$ is minimal in the straightening relation; in particular $\nu(\sigma_{h-1}\pi) \vgeq \nu(\pi_{j,1}\pi_{j,2})$.

Now we apply \myref{Lemma}{lemma_impossibleCase} (to $\pi = \pi$, $\pi_1 = \pi_{j,1}$ and $\pi_2 = \pi_{j,2}$ with notation as in that Lemma) and find that there exists $\tau$ such that $\pi(\tau) < \pi_{j,2}(\tau)$ with $k=\ell(\tau)\geq h$ and such that $\pi(\sigma) = \pi_{j,2}(\sigma)$ for each $\sigma > \tau$. By $\pi_{j,2}(\tau) > \pi(\tau)$ we have
\[
\pi_{j,2}(\tau)\geq \pi(\tau) + \frac{1}{M_\tau} ,
\]
hence $\nu_{0,\calC}(\pi_{j,2})\geq\nu_{0,\calC}(\pi) + \sigma'_h + \bar{v}$, for some $\bar{v}\in V_{h-1}$. So we find
\[
\sigma_{h-1} + \nu(\pi) = \nu(\sigma_{h-1}\pi) \vgeq \nu(\pi_{j,1}\pi_{j,2}) \vgeq \nu(\pi_{j,2}) \vgeq \nu_{0,\calC}(\pi_{j,2}) \vgeq \nu_{0,\calC}(\pi) + \sigma'_h + \bar{v},
\]
thus, setting $v'=-\sigma'_{h-1} + \bar{v}\in V_{h-1}$, we have $\nu(\pi) \vgeq \nu_{0,\calC}(\pi) + \sigma'_h + v'$; our claim is proved.
\end{proof}

\subsection{Quasi-valuations for LS algebras}\label{subsection_quasiValuation}

In this subsection we introduce the main object of the paper: the quasi-valuation defined as the minimum of positive valuations $\nu_\calC$ with prescribed values along the chain $\calC$, with $\calC$ running over the set of all maximal chains in $S$. This quasi-valuation will have values in $\Q^S$, so we need a total order on $\Q^S$: we consider a total order on $S$ refining the given partial order and let $\vleq$ be the reverse lexicographic order on $\Q^S$ induced by this order.

Suppose that, for each maximal chain $\calC$ in $S$, we are given a \emph{positive} valuation $\nu_\calC:A\setminus\{0\} \longrightarrow \Q^\calC\subseteq \Q^S$ such that $\nu_\calC(\sigma) = \sigma$ for each $\sigma\in\calC$. Note that if $A$ is of flag type then, by \myref{Theorem}{theorem_positivityForFlagType}, such positive valuations exist.

Also, by \myref{Theorem}{theorem_estimate}, if $A$ has an effective system of weights (e.g. it is of flag type) then for any LS path $\pi$ we have $\nu_{\calC}(\pi)\geq \nu_{0,\calC}(\pi) + \sigma'_{h_\calC(\pi)} + v'$, for some $v'\in V_{h_\calC(\pi)-1,\calC}$; in particular we have $\nu_\calC(\pi) \geq \nu_{0,\calC}(\pi)$ and $\nu_\calC(\pi)=\nu_{0,\calC}(\pi) = \pi$ if and only if $\supp\pi\subseteq\calC$.

Finally note that if we restrict $\vleq$ to $\Q^\calC$, with $\calC$ a chain in $S$, we get the order we have used in the previous subsections; so the symbol $\vleq$ is not ambiguous.

We define the main object of interest of this Subsection:
\[
\begin{array}{rcl}
A\setminus\{0\} & \stackrel{\nu}{\longrightarrow} & \Q^S\\
x  & \longmapsto & \min_{\calC}\nu_{\calC}(x),\\
\end{array}
\]
where the minimum runs over all maximal chains in $S$ and is taken with respect to the total order $\vleq$ in $\Q^S$. Being the minimum of a finite number of valuations, the map $\nu$ is a quasi-valuation by \myref{Proposition}{proposition_quasiValuationAsMin}.

In general, the quasi-valuation $\nu$ does depend on the fixed total order on $S$ used to define $\vleq$. However, in this subsection we show that its values on the LS paths is independent of this choice.

In the following Lemma we compare valuations associated to different chains.
\begin{lemma}\label{lemma_compareTwoChains} Let $\pi\in\LS$ be an LS path with support contained in the maximal chain $\calC$ and \emph{not} contained in the maximal chain $\calD$. Then $\nu_\calD(\pi) \vgneq \nu_\calC(\pi) = \pi$.
\end{lemma}
\begin{proof} Let $h=h_\calD(\pi)$, let $\sigma$ be the $h$--th element of the chain $\calD$ and let $\sigma' = \sigma / M_\sigma$. By \myref{Theorem}{theorem_estimate}, $\nu_\calC(\pi) = \nu_{0,\calC}(\pi) = \pi$ while $\nu_\calD(\pi)=\nu_{0,\calD}(\pi) + \sigma' + v'$ for some $v'\in V_{h-1,\calD}$. By the definition of $h$ we have $\nu_{0,\calC}(\pi_{|S_{>\sigma}})=\nu_{0,\calD}(\pi_{|S_{>\sigma}}$). So we deduce
\[
\begin{array}{rcl}
\nu_\calD(\pi) & \vgeq & \nu_{0,\calD}(\pi) + \sigma' + v'\\
 & = & \nu_{0,\calD}(\pi_{|S_{>\sigma}}) + (\pi(\sigma) + 1/M_\sigma)\sigma + (\nu_{0,\calD}(\pi_{|S_{<\sigma}}) + v')\\
 & \vgneq & \nu_{0,\calC}(\pi_{|S_{>\sigma}}) + \pi(\sigma)\sigma' + \nu_{0,\calC}(\pi_{|S_{<\sigma}}), 
\end{array}
\]
where we have used that $\sigma' > -(\nu_{0,\calD}(\pi_{|S_{<\sigma}}) + v') + \nu_{0,\calC}(\pi_{|S_{<\sigma}})$.
\end{proof}

We are now in a position to prove the following fundamental result.

\begin{theorem}\label{theorem_quasiValuation} The quasi-valuation $\nu$ has the following properties:
\begin{itemize}
\item[(a)] $\nu(\pi)=\pi$ for any LS path $\pi$ and
\[
\nu\Big(\sum c_\pi\pi\Big) = \min_\vleq\Big\{\pi\,|\,c_\pi\neq 0\Big\},
\]
in particular, $\nu$ depends only on the fixed total order on $S$ used to define $\vleq$,
\item[(b)] for any monomial in LS paths $\pi_1\pi_2\cdots\pi_r$ and any maximal chain $\calC$, we have $\nu(\pi_1\pi_2\cdots\pi_r) \vleq \nu_\calC(\pi_1\pi_2\cdots\pi_r)$, where the equality holds when $\supp\pi_h\subseteq\calC$ for $h=1,2,\ldots,r$,
\item[(c)] if the LS paths $\pi_1,\pi_2,\ldots,\pi_r$ have comparable supports then
\[
\nu(\pi_1\pi_2\cdots\pi_r) = \pi_1 + \pi_2 + \cdots + \pi_r,
\]
\item[(d)] if the LS paths $\pi_1,\pi_2,\ldots,\pi_r$ have non-comparable supports then
\[
\nu(\pi_1\pi_2\cdots\pi_r) \vgneq \pi_1 + \pi_2 + \cdots + \pi_r.
\]
\end{itemize}
\end{theorem}
\begin{proof}\begin{itemize}
\item[(a)] Let $\calC$ be a maximal chain in $S$. By \myref{Lemma}{lemma_compareTwoChains}, if $\calC$ contains $\supp\pi$, then $\nu_\calC(\pi) = \nu_{0,\calC}(\pi) = \pi$, otherwise $\nu_\calC(\pi) \vgneq \nu_{0,\calC}(\pi)$. Hence $\nu(\pi)=\min_{\calC}\nu_\calC(\pi)=\pi$ since there exists at least one maximal chain containing the support of $\pi$.

For the second claim note that $\nu$ is injective on LS paths by what just proved. Thus, by \myref{Lemma}{lemma_quasiValuationMin},
\[
\nu(\sum c_\pi\pi)=\min_\vleq\{\pi\,|\,c_\pi\neq 0\}
\]
since all the values $\nu(\pi)$'s are different.
\item[(b)] The inequality is clear by the definition of $\nu$. On the other hand, using again \myref{Lemma}{lemma_compareTwoChains}, if $\supp\pi_h\subseteq\calC$ for $h=1,\ldots,r$, then $\nu_\calD(\pi_1\cdots\pi_r) \vgeq \nu_\calD(\pi_1) + \cdots + \nu_\calD(\pi_r) \vgeq \nu_\calC(\pi_1) + \cdots + \nu_\calC(\pi_r) = \pi_1 + \cdots + \pi_r$ for any maximal chain $\calD$, thus $\nu(\pi_1\cdots\pi_r) \vleq \nu_\calC(\pi_1\cdots\pi_r)$.
\item[(c)] This follows at once by applying (b) to a maximal chain containing the supports of all $\pi_1,\ldots,\pi_r$.
\item[(d)] For any maximal chain $\calC$ there exists $1\leq h\leq r$ such that $\supp\pi_h\not\subseteq\calC$. By Theorem \myref{Theorem}{theorem_estimate} $\nu_\calC(\pi_j) \vgeq \pi_j$ for any $j=1,\ldots,r$ and $\nu_\calC(\pi_h)  \vgneq \pi_h$, hence $\nu_\calC(\pi_1\cdots\pi_r) = \nu_\calC(\pi_1)+\cdots+\nu_\calC(\pi_r) \vgneq \pi_1+\cdots+\pi_r$.

Being $\nu(\pi_1\cdots\pi_r)$ the minimum over the finite set of maximal chains of $S$, we have also $\nu(\pi_1\cdots\pi_r) \vgneq \pi_1+\cdots+\pi_r$.
\end{itemize}
\end{proof}

\section{Consequences}\label{section_consequences}

In the previous section we have constructed a quasi-valuation $\nu$ as the minimum of positive valuations along maximal chains. This allows an interpretation of the LS algebra structure, i.e. of the standard monomial theory for $A$, in terms of the quasi-valuation. First of all the set of LS paths itself gets such an interpretation.
\begin{corollary}\label{corollary_LSPathsIsImage}
The image of $\nu$ in $\Q^S$ is the set of LS paths.
\end{corollary}
\begin{proof} This follows at once by (a) of \myref{Theorem}{theorem_quasiValuation}.
\end{proof}

The degeneration of an LS algebra to the discrete algebra (see \cite{chiriviLS}) is just an instance of the general result of degeneration arising from the graded algebra associated to the quasi-valuation.

Indeed we can recover the discrete LS algebra over $(S,\leq,\bond)$ as the graded algebra associated to the quasi-valuation. Indeed, for $\pi\in\LS$ let $A_{\vgeq\pi} = \{x\in A\,|\,x = 0\textrm{ or }x \neq 0\textrm{ and }\nu(x) \vgeq \pi\}$ and $A_{\vgneq\pi} = \{x\in A\,|\,x = 0\textrm{ or }x\neq 0\textrm{ and }\nu(x) \vgneq \pi\}$. Since $\LS\subseteq\Q^S_{\geq0}$ the $\K$--vector space
\[
\gr_\nu(A) = \bigoplus_{\pi\in\LS} A_{\vgeq\pi} / A_{\vgneq\pi}
\]
is an algebra graded by $\LS$, this is called the \emph{degenerate algebra} of $A$ with respect of $\nu$.
\begin{corollary}\label{corollary_degenerateAlgebra}
The degenerate algebra $\gr_\nu(A)$ is isomorphic to the discrete LS algebra over $(S,\leq,\bond)$.
\end{corollary}
\begin{proof} Since $\LS_1$ generates $A$ as a $\K$--algebra, its image in $\gr_\nu(A)$ is yet a set of generators; we denote by $\pi$ also the image in $\gr_\nu(A)$ of the LS path $\pi$.

By \myref{Theorem}{theorem_quasiValuation} the monomial $\pi_1\pi_2$, with $\pi_1,\pi_2\in\LS_1$, is $0$ in $\gr_\nu(A)$ if $\pi_1,\pi_2$ have non-comparable supports. On the other hand let $\pi_1\pi_2$ be a non-standard monomial with $\pi_1,\pi_2$ having comparable supports, let $\pi_{0,1}\pi_{0,2}$ be its canonical forms and let
\[
\pi_1\pi_2 = \pi_{0,1}\pi_{0,2} + \sum_{j}c_j\pi_{j,1}\pi_{j,2}
\]
be its straightening relation in $A$.

By \myref{Theorem}{theorem_quasiValuation} we know that $\nu(\pi_1\pi_2) = \min_j\nu(\pi_{j,1}\pi_{j,2})$ where the minimum is taken with respect to $\lsleq$. But $\nu(\pi_{j,1}\pi_{j,2}) = \pi_{j,1} + \pi_{j,2}$ since $\pi_{j,1}$ and $\pi_{j,2}$ have comparable supports for any $j$. On the other hand, $\pi_{0,1} + \pi_{0,2} \lslneq \pi_{j,1} + \pi_{j,2}$ for all $j > 0$. Since $\lsleq$ is stronger than $\vleq$ we conclude $\nu(\pi_1\pi_2) = \nu(\pi_{0,1}\pi_{0,2})$ so that $\pi_1\pi_2 = \pi_{0,1}\pi_{0,2}$ in $\gr_\nu(A)$.
\end{proof}

Another direct consequence is that we can explicitly compute the Newton-Okounkov body of the quasi-valuation, it is just the order complex of the poset. This follows immediately by \myref{Corollary}{corollary_LSPathsIsImage}.

\begin{corollary}\label{corollary_NOquasiValuation}
The Newton-Okounkov body $\NO_\nu(A)$ of $A$ with respect to the quasi-valuation $\nu$ is the order complex $\Delta(S)$ of $(S,\leq)$. In particular, the collection of levels of $\NO_\nu(A)$ is an integral structure $\{\Delta_r(S)\,|\,r=1,2,\ldots\}$ of $\Delta(S)$.
\end{corollary}

\section{Application to Schubert Varieties}\label{section_application}

\subsection{Schubert varieties and LS algebras}

Let $\lieg$ be a semisimple Lie algebra (or a symmetrizable Kac-Moody algebra) over an algebraically closed field $\K$ of characteristic $0$ and let $\lieh\subseteq\lieb\subseteq\lieg$ be a Cartan subalgebra and a Borel subalgebra, respectively. Let $\Lambda$ be the weight lattice of $\lieg$ and let $\Lambda^+$ be the monoid of dominant weights corresponding to $\lieb$. Let $(\cdot,\cdot)$ be the Killing form on $\Lambda\otimes_\Z\R$ and $\beta^{\vee}=2\beta/(\beta,\beta)$ the coroot corresponding to the root $\beta$. Fix a dominant weight $\lambda$ and let $W_\lambda$ be its stabilizer in the Weyl group $W$ of $\lieg$; denote by $\leq$ the restriction of Bruhat order to the set of minimal representatives $W^\lambda\subseteq W$ of the quotient $W/W_\lambda$.

Let $G$ be the simply connected semisimple algebraic group with Lie algebra $\lieg$, let $B\subseteq G$ be the Borel subgroup with Lie algebra $\lieb$ and let $P\supseteq B$ be the parabolic subgroup stabilizing the dominant weight $\lambda$. The Schubert variety $X(\tau)$, for $\tau\in W/W_\lambda$ is the closure of the $B$--orbit $B\tau P/P$ in the partial flag variety $G/P$. Note that, if $\lieg$ is of finite type, then $X(w_{0,\lambda})=G/P$ where $w_{0,\lambda}$ is the longest element in $W/W_\lambda$. In the sequel we fix a $\tau\in W/W_\lambda$.

Let $\calL_\lambda$ be the line bundle $G\times_P\K_{-\lambda}$ on $G/P$; recall that the space $H^0(G/P,\calL_\lambda)$ of global sections is the $G$--module $V(\lambda)^*$ dual of the module $V(\lambda)$ of highest weight $\lambda$. Denote also the restriction of $\calL_\lambda$ to $X(\tau)$ by $\calL_\lambda$ and consider the graded $\K$--algebra
\[
A = A_0\oplus A_1\oplus\cdots,\quad A_n \doteq H^0(X(\tau),\calL_{n\lambda}).
\]

Let $W^\lambda_\tau$ be the set of all $\sigma\leq\tau$ in $W^\lambda$; this is a poset with the restriction of the Bruhat order. Moreover $W^\lambda_\tau$ has $e$ as the unique minimal element, $\tau$ as the unique maximal element and it is a graded poset with length $N=\dim X(\tau)$.

If $\sigma'$ covers $\sigma$ in $W^\lambda_\tau$ then there exists a positive root $\beta$ such that $\sigma'=s_\beta\sigma$, with $s_\beta\in W$ the symmetry with respect to $\beta$. We define $\bond_\lambda(\sigma,\sigma')=(\sigma(\lambda),\beta^{\vee})$; this is a set of bonds for the poset $W^\lambda_\tau$; the $\gcd$ condition on chains follows by Proposition~2.7 in \cite{dehy}.

In a series of papers \cite{L1}, \cite{L2}, \cite{L3}, \cite{L4} the third author introduced the language of LS paths for symmetrizable Kac-Moody algebras. In his notation an LS path is a pair of sequences $(\sigma_1>\cdots>\sigma_r\,;\,0=a_0<a_1<\cdots<a_r=1)$ with $\sigma_1,\ldots,\sigma_r\in W^\lambda$ and rational numbers $a_0,\ldots,a_r$ with certain integral conditions. As proved in details in \cite{chiriviLS} the map
\[
(\sigma_1>\cdots>\sigma_r\,;\,0=a_0<a_1<\cdots<a_r=1) \longmapsto \sum_{j=1}^r (a_j - a_{j-1})\sigma_j
\] 
is a bijection between Littelmann's LS paths and LS paths of degree $1$ for $(W^\lambda_\tau,\leq,\bond_\lambda)$ as defined in \myref{Section}{section_LSAlgebra} (see also \cite{CFL2} for this proof). Moreover in \cite{L5}, a section $p_\pi$ in $H^0(X(\tau),\calL_{\lambda})=A_1$ is defined for each LS path $\pi\in\LS_1(W^\lambda_\tau)$. 

In \cite{L5} and in \cite{LLM} certain relations for the section $p_\pi$'s are proved. As explained in \cite{chiriviLS}, in our setting, this may be stated as: $A$ is an LS algebra over the poset with bonds $(W^\lambda_\tau,\leq,\bond_\lambda)$.

Moreover the map
\[
\LS\,\ni\,\pi\,\longmapsto\,\sum_{\sigma\in W^\lambda_\tau}\pi(\sigma)\sigma(\lambda)\,\in\,\Lambda
\]
is an effective system of weights for $A$: the straightening relations of $A$ are homogeneous for $\Lambda$ since the torus $T\subseteq B$ corresponding to $\lieh$ acts on $X(\tau)$ and the system is effective since the extremal LS paths of degree $n$ correspond to extremal weight vectors in $H^0(G/P,\calL_{n\lambda}) \simeq V(n\lambda)^*$ which have multiplicity $1$.

Let $v_\lambda$ be a highest weight vector in $V(\lambda)$, consider the embedding $G/P\ni gP\longmapsto [g\cdot v_\lambda]\in \Pj(V(\lambda))$ and its restriction $X(\tau)\longrightarrow\Pj(V(\lambda)_\tau)$ where $V(\lambda)_\tau$ is the vector subspace of $V(\lambda)$ spanned by $X(\tau)$.

This embedding of $X(\tau)$ identifies the line bundle $\calL_\lambda$ with $\mathcal{O}(1)$ on $\Pj(V(\lambda)_\tau)$. In particular the component $V(\lambda)_\tau^*$ of degree $1$ of the coordinate ring $\K[\hat{X}(\tau)]$ of the cone $\hat{X}(\tau)$ over $X(\tau)$ in $V(\lambda)_\tau$ is isomorphic to $A_1$ as a $B$--module.

Now $\K[\hat{X}(\tau)]$ is clearly generated by $V(\lambda)_\tau^*$ and $A$ is generated by $A_1$ since it is an LS algebra. In conclusion we have an isomorphism $\K[\hat{X}(\tau)]\simeq A$ as $B$--modules. In the sequel we will freely identifies these two algebras.

As proved in \cite{L5}, for $\sigma\in W^\lambda_\tau$, the ideal $I_\sigma$ of $A$ generated by the LS paths $\pi$ such that $\max\supp\pi\not\leq\sigma$ is a defining ideal for (the cone over) $X(\sigma)$ in $X(\tau)$. Hence $A$ is of flag type since the quotient $A_\sigma = A / I_\sigma$ is thus the coordinate ring of the cone over the embedding of the Schubert variety $X(\sigma)\subseteq X(\tau)\hookrightarrow\Pj(V_\lambda)$ and it is a domain being $X(\sigma)$ irreducible.

We summarize this discussion in the following theorem.
\begin{theorem}\label{theorem_SchubertLS}
The coordinate ring $A$ of the cone over the embedding $X(\tau)\hookrightarrow\Pj(V_\lambda)$ of the Schubert variety $X(\tau)$ is an LS algebra of flag type over the poset with bonds $(W^\lambda_\tau,\leq,\bond_\lambda)$ with respect to the injection
\[
\LS\,\ni\,\pi\,\longmapsto\,p_\pi\,\in\, A.
\]
\end{theorem}

\myref{Theorem}{theorem_positivityForFlagType} implies that for each maximal chain $\calC$ in $W^\lambda_\tau$ there exists a positive valuation $\nu_\calC:A\setminus\{0\}\longrightarrow\Q^\calC\subseteq\Q^{W^\lambda_\tau}$ such that $\nu_\calC(p_\tau) = \tau$ for each $\tau\in\calC$. So \myref{Theorem}{theorem_estimate} can be applied to the quasi-valuation $\nu = \min_\calC\nu_\calC$. In particular we get an interpretation of the LS paths, a semi-toric degeneration and the description of the Newton-Okoukov body.
\begin{theorem}\label{theorem_quasivaluationSchubert}
\begin{itemize}
	\item[(i)] The image of $\nu$ is the set of LS paths $\LS$ for the poset with bonds $(W^\lambda_\tau,\leq,\bond)$.
	\item[(ii)] The quasi-valuation $\nu$ defines a flat degeneration of the Schubert variety $X(\tau)$ to a union of irreducible normal toric varieties, one for each maximal chain in $W^\lambda_\tau$.
	\item[(iii)] The Newton-Okounkov body $\NO_\nu(A)$ of $A$ with respect to the quasi-valuation $\nu$ is equal to the order complex $\Delta$ of the set of minimal representatives $W^\lambda_\tau$ ordered by the Bruhat order. The chains of $W^\lambda_\tau$ defines a triangulation of $\NO_\nu(A)$ and the collection of levels of $\NO_\nu(A)$ is an integral structure of $\Delta$.
\end{itemize}
\end{theorem}
\begin{proof} The first claim is \myref{Corollary}{corollary_LSPathsIsImage}. The second claim is a consequence of \myref{Corollary}{corollary_degenerateAlgebra} since the discrete algebra is the coordinate ring of the cone over a union of irreducible normal toric varieties, one for each maximal chain, as proved in \cite{chiriviLS}. The last claim is \myref{Corollary}{corollary_NOquasiValuation}. 
\end{proof}

\subsection{The case of Grassmannians} Let $\lambda=\omega_k$ be a fundamental weight for the Lie algebra $\lieSL_{\ell+1}(\K)$. Then $G/P$ is the Grassmannian of $k$--spaces in $\K^{\ell+1}$. This case is studied in \cite{fangLittelmann}. For such $\lambda$ all bonds are equal to $1$, i.e. we are in the Hodge algebra case, the poset $W^\lambda$ is a distributive lattice. LS paths of degree $1$ are just elements of $W^\lambda$. Moreover, the weight spaces of $V_\lambda$ are all $1$--dimensional; so the path basis is uniquely determined up to non-zero constant factors, hence it is simply the basis of Pl\"ucker coordinates. In particular the path basis coincides with the crystal/canonical basis (see \cite{lusztig} for canonical basis).

This example may be verbatim generalized to any minuscule weight $\lambda$.

As mentioned in the Introduction, the second and third author in \cite{fangLittelmann} defined certain valuations $\nu_\calC$, with $\calC$ a maximal chain in $S$, and the quasi-valuation $\nu$ as the minimum of the $\nu_\calC$'s. This quasi-valuation as the same properties of the quasi-valuation defined in the present paper. However their construction is different from the one in this paper; that in \cite{fangLittelmann} relies in an essential way on the structure of distributive lattice of $W^\lambda$ for minuscule $\lambda$. Moreover they use only partially the structure of Hodge algebra, in particular they do not need the order property of the straightening relations but add some hypothesis about certain standard monomials appearing in the straightening relations.

\subsection{Compatibility with valuation from the Seshadri stratification}
In \cite{CFL} the notion of Seshadri stratification for a projective variety is introduced. We briefly recall this notion in our context of Schubert varieties; see \cite{CFL} for generalities about Seshadri stratifications and also for the details about this application.

We keep denoting by $\tau$ a fixed element in $W^\lambda$. For $\sigma\in W^\lambda_\tau$ let $p_\sigma$ be the \emph{extremal function} associated to the extremal path $\sigma$. As proved in \cite{CFL}, the collection of Schubert subvarieties $X(\sigma)\subseteq\Pj(V_\lambda)$, and extremal functions $p_\sigma$, $\sigma\in W^\lambda_\tau$, define a Seshadri stratification of $X(\tau)$ since the following properties hold:
\begin{itemize}
	\item[(S1)] for each $\sigma$, the subvariety $X(\sigma)$ is non-singular in codimension one; for each covering relation $\sigma < \sigma'$, $X(\sigma)$ is a codimension one subvariety of $X(\sigma')$,
	\item[(S2)] if $\sigma\not\leq\sigma'$ then $p_\sigma$ vanishes on $X(\sigma')$,
	\item[(S3)] the set of zeros of $p_\sigma$ on $X(\sigma)$ is set-theoretically the union of the $X(\sigma')$'s with $\sigma'$ covered by $\sigma$.
\end{itemize}

In the same paper, a valuation $\nu_\calC:A\setminus\{0\}\longrightarrow\Q^\calC\subseteq\Q^{W^\lambda_\tau}$ is defined for each maximal chain $\calC$ in $W^\lambda_\tau$. These valuations are positive and have the property: $\nu_\calC(p_\sigma) = \sigma$ for each $\sigma\in\calC$. So the results of the present paper, in particular \myref{Theorem}{theorem_estimate}, can be applied to these valuations.

Moreover, in \cite{CFL} a quasi-valuation $\nu$ is defined as the minimum of the valuations $\nu_\calC$, $\calC$ a maximal chain, using a total order $\leq^t$ on $W^\lambda_\tau$ refining the Bruhat order $\leq$ and such that $\sigma<^t\sigma'$ if $\ell(\sigma) < \ell(\sigma')$. Of course we may use this order also in the construction of the quasi-valuation as in \myref{Subsection}{subsection_quasiValuation}; so \myref{Theorem}{theorem_quasiValuation} and \myref{Corollary}{corollary_LSPathsIsImage} can be applied to $\nu$. We get the following very interesting geometric interpretation of the combinatorial notion of LS path.
\begin{theorem}
For each LS path on $W^\lambda_\tau$ we have $\nu(p_\pi) = \pi$. In particular, the LS paths on $W^\lambda_\tau$ are the indexing set of the non-zero leaves of a quasi-valuation defined as the minimum of certain valuations $\nu_\calC$, $\calC$ a maximal chain in $W^\lambda_\tau$.
\end{theorem}

So the LS paths encode the vanishing orders of certain regular functions on the web of the Schubert subvarieties of $X(\tau)$. This theorem is proved also in \cite{CFL2} using the construction of the section $p_\pi$ via quantum groups in a more direct manner.

\end{document}